# Title: Evolutionary Understanding of the Conditions Leading to Estimation of Behavioral Properties through System Dynamics

**Authors:** Chulwook Park

**Affiliations:**[1]International Institute for Applied Systems Analysis (IIASA), A-2361 Laxenburg, Austria. *To whom correspondence should be addressed. E-mail: parkc@iiasa.ac.at

**Abstract:** One of the basic frameworks in science views behavioral products as a process within a dynamic system. The mechanism might be seen as a representation of many instances of centralized control in real time. Many real systems, however, exhibit autonomy by denying statically treated mechanisms. This study addresses the issues related to the identification of dynamic systems and suggests how determining the basic principles of a collective structure may be key to understanding complex behavioral processes. A fundamental model is derived to assess the advantages of this perspective using a basic methodology. The connection between perspective and technique demonstrates certain aspects within their actual context, while also clearly including the framework of actual dynamic system identification.

**One Sentence Summary:** This study addresses the issues related to the identification of dynamic systems and suggests how determining the basic principles of a collective structure may be key to understanding complex behavioral processes. Distinct approaches (theoretical, modelling, and experimental) were used in an effort to understand how to recognize the basic feature through the dynamics of the system. The condition leading to the behavioral property under these different perspectives suggests that it is not reliant on an individual lodged in the corner of the system. Rather, if the fundamental condition is in place, the complexity of the system falls into the same rules that are estimating the pattern.

**Main Text:** System dynamics can be as a tool to address many applications in a very broad sense. It applies to systems in many different disciplines, such as control, communication, tropisms, and even systems of systems. Given the very widespread use of the term "system," the important questions are: What does the term "system" actually mean? and What makes a system "dynamic"? The present study does not claim to produce a rigid and rigorous definition; instead, it attempts to recognize the meaning intuitively through a general discussion.

A system is an abstract or vague entity. It would not be defined through strict codification (*1*). A well-formulated notion of a system could, however, be developed using a non-system definition. A non-system can be represented by a set of isolated entities that do not interact with each other, or a collection of entities whose relationships have no implications for the properties and behaviors of the entities (*2*). Given how vague word definitions can be, the notions that we use describe a wide variety of things (*3*). Let us thus specify how these notions may be understood in the light of an account of a dynamic system.

First, a system is a collection of individual elements: it is not usual to talk about a system having only one component (*4*). It is then likely that the elements comprising the system exhibit non-trivial interactions, for instance, coherence in terms of moving or working together [$\chi \rightarrow \vec{f} \rightarrow f(x_i)$]. Here, $x$ denotes an element composed of many individuals; $f$ is a function toward a certain direction ($\rightarrow$), such as a goal, and $f(x_i)$ denotes coherence, which has certain individuals ($x_n$), but not all individuals, underlying a very well-defined goal boundary ($f$). The second important—and key—assumption of the system is that the well-defined individuals must interact not just within one another but also with everything outside the system (*5*), which is called the environment $x_t \rightarrow f(x_i) \rightarrow y_t$. In this equation, there are two kinds of interactions that $x_i$ can have with the outside environment. The outside exerts some impact ($x$) on them and is influenced ($y$) by them at a certain time ($t$). Thus, intuitively, a system can be defined as many interactive individual elements embedded within a certain coherent behavior (*6, 7*).

### *Part 1: Nonrepresentational perspective of phenomena (theoretical basis)*
There is a well-defined assumption called the problem of impoverished entailment (*8*), which is the minimal starting point for understanding any system of interest at any level of interest. "X is about Y" is true only if "X entails Y" and "Y entails X" are true. This is a loop of entailment; specificity of X to its source of Y means that X entails the source by which X is entailed. Diagrammatically, an entailment can be expressed as X→Y. The primary property of an entailment is that it propagates "truth" hereditarily—Y inherits the truth of X (*9*). Thus, the loop of entailment can be shown diagrammatically as X→Y →X, and the diagram can be read as meaning that truth propagates hereditarily in both directions. Such an explanation focuses on assembling all parts of elements [$C(s)$] along with the other "things" that can influence them or be influenced by them [$E(s)$], even the internal structure of the system [$S(s)$] (*10*).

$$m(s) = < C(s), E(s), S(s) > \tag{1}$$

In this framework, when $m$ is predicated as a mass of a particular object, it identifies a substantial property intrinsic to the object that is identical whenever and wherever the object is observed. $m$ has other relational properties, however, that engage actual things. For example, when the object is grasped and brought down hard and repetitively on another object, then $m$ can be a "hammer" with respect to the object. Conceivably, there are many relational properties that the object may have by virtue of its relationship either with other objects or with perceivers-actors, but these are, at best, indefinite properties until a particular spatio-temporal relationship is effected (*11*). When this occurs, one of the many potential relational properties of the object is actualized.

However, investigating the dynamics in these cases typically requires a significant number of elements and obviously includes multiple components that must be managed. According to researchers (*12*), organisms have too many degrees of freedom (humankind = $10^{19}$). Moreover, if the parts are considered to be very strongly defined by their connections and to function within the context, great complexity can to be observed (*13*). When even a single cell's behavior is being considered, its tendencies are not certain to predict all of these dynamics (*14, 15*), and this poses deep concerns over how to treat the cell relevantly as well as questions about how these behaviors form from numerous factors.

Instead of positively interpreting the system's dynamics in terms of every encoded equation, the key idea is to describe why numerous modes of emergent phenomena underlying local-level interaction have to be governed by simple rules (*16*). Where elements are given a simple set of rules that govern their behavior and allow them to interact to determine what patterns emerge over time, it has been shown that an agent's behavior with respect to unpredictable phenomena can arise even with elementary governing rules. This fresh perspective on refocusing a system's dynamics is helping to bridge traditional biases and has been stimulating scientists to distil out simple principles, such that a better understanding of a dynamic system (which is not always very complicated) may be gained. The following section constitutes syntheses from investigations of various behavioral phenomena in terms of fundamental mechanisms that can be depicted on a specific screen, and that obtain simplicity from complexity (*13*).

### Part 2: Individual behavior in social dynamics (Model-based)

To explore the rule of thumb (*17*) from a broader perspective, computer simulation was performed first on the basics of a spatially explicit model of mobile agents in continuous space to determine the basic regulatory principles involved in the way they conceptualize their environment (*18*). In other words, the model implemented shows the potential to infer how simple individual rules can lead to consistent group behavior, and also how slight changes in those mechanisms can have a dramatic impact on an individual's behavioral patterns (*19*). Although an analysis of simple implication is an apparent first step in providing proof of concept, this individual-based simulation has become significant enough to be tested in a broader range of applications for evolutionary dynamics for the following reasons. First, as the agents represent individuals that have occurred from the bottom up, the actual state of their behavior tends to be more informative (agent-based modelling). Next, as the main point of this implementation is a description that deals with state per time as the critical factor in its allocation of neighbors, the number of neighbors placed on the position is based on the moves scheduled for a given moment (context-varying cultural evolution) (*20*) (see Supplement 1.2 for more detail).

### Result

The model provides a natural description of a pattern of behavior, and allows to understand a realistic adaptation incorporating behavioral algorithms with social dynamics. The mechanisms characterized in the agenda shows that the model has the capacity to produce three types of factor. First, the behavioral pattern is the result of the applied individual components; this comes not only from the initial conditions of the autonomous agents but also from the fact that they are interconnected. Second, the range of different combinations of the internal and the external state plays a part in the rapid propagation in the system (Fig. 1.1). Third, however, when the interconnected relation between the internal trait underlying its external trait is applied, the widespread heterogeneousness of the mechanisms can abstract the repertoire of displayed behaviors (Fig. 1.2) (see Supplement 1.3 for more detail).

The primary feature of the agents' interactions is heterogeneous in this abstract setting. As the topology of the interaction traits can lead to significant deviations from the predicted pattern of behavior, it may generate various effects that mimic the behavior of real individuals in social dynamics. At the points where these individuals interact, sensible decisions occur, in line with an empirical study showing that individuals would learn how to keep relative velocity as a key factor for homogeneousness (*21*).

### Discussion

The simulations show that different interconnection structures have an effect on which strategies perform better: a relationship referred to as ecological rationality. The results may propose that relative velocity can be the pure candidate and works. According to researchers (*21*), if an object (i.e., ball) is already high in the air and travelling directly in line with the individual (i.e., player), the individual might utilize some simple heuristics. Namely, the individual fixes his gaze on the object, starts running, and adjusts his velocity to ensure that the angle of the ball above the horizon appears constant (*22*). The prediction is not that the individual runs to a pre-computed landing spot and waits for the object, but that he is modified to keep the image of the object moving at a constant velocity. It is also possible that individuals do not compute ($\vec{v}$) at all in this model but would reduce a maintained value of $d^2(\vec{v})/dt^2$ in a systematic way. As $\vec{v}$ increased, they would keep $d^2(\vec{v})/dt^2$ at zero [$d^2(\vec{v})/dt^2 = $ constant] (*23*). This is related to what information can be derived as a strategy in the system and how that information can be best obtained through the dynamics.

This result also holds for a progressed step toward such a relationship, as explicitly shown by the dramatic change when it comes to the point where [$d^2(\vec{v})/dt^2 = $ constant] reaches [$d^2(\vec{v})/dt^2 \neq $ constant]. Let us suggest that actualized observation of the pattern corresponds not simply to the object's velocity but, instead, to changes in the velocity between individuals at a certain point. In other words, the actual displacement ($S$) estimated is given by ($S_{l^1} = S_l + \alpha d^2$). Here, the observed new displacement ($S_{l^1}$) is equal to the displacement across the individuals ($l^1$) plus its relative velocity ($\alpha$) multiplied by the distance squared ($d^2$) from the neighbor. This indicates that the further away an individual is from the neighbor (role model or group heading), when measured at a certain point, the more difficult it is for the individual to follow the neighbor because of the greater acceleration involved. On this assumption, the results might be able to suggest a strategy such as evolvable traits or payoff functions to guide the evolution of these heuristics through social learning (Fig. 2).

***Part 3: Elementary coordination in circadian rhythm (experimental-based)***

Let us look at how, each day, we carry out very basic actions with the primary goal of obtaining simplicity from complexity. We can move segments of our body (fingers, arms, legs, arms and legs, head and arms, and so on). Segments and limbs of other organisms, too, such as horses, dogs, and even fish, are all involved in locomotion (*24*). The phenomenon of how animals typically use their limbs in a number of distinct periodic modes is familiar to us; moreover, many of these modes possess some degree of symmetry (*25*). This is a very fundamental point to understand, given that animals, as they travel, continue to use these modes, at some points at speed.

How, then, ought we to understand this particular ability? One useful strategy looks for cycles at all time scales and aims to show how interacting cyclic processes account for the emergence of new entities (*26, 27*), many of which are similarly cyclic (*28*). The central idea is that the Earth's cycles—geophysical, hydrological, meteorological, geochemical, and biochemical—have interacted to create self-replicating living systems that abide by particular cyclicities (*29*). This assumption has led us to enquire whether something akin to attunement to the environmental 24-hour day/night cycles (*30*) may be apparent in an experimental setting of bimanual coordination, a context that has been used to examine self-organization in biological systems (*31*).

The present study used two main ways of determining these characteristics and discovering if approximations under certain conditions serve these self-potentials (see Supplement 2.1 for more detail). The first involves an increase in the capability to self-generate forces along the lines of the roles of the fundamental dimensions of environments (temperature embedding in light-dark cycles). To achieve this, the experimental setting asks, "Is our system influenced by an ecological feature?" by embedding a bimanual coordination task in an ordinary 24-hour day–night cycle (5:00, 12:00, 17:00, and 24:00). The second is tied to observing the availability of an internally based source (coordination) or sources of force (stability and entropy) within dynamical boundaries in systematic ways. The setting asks, "How does our system adapt to regular or irregular thermal structures?" by embedding the comparison of normal and abnormal day–night circadian temperature effects at dawn (5 a.m., approximately when core temperature reaches its minimum) and dusk (5 p.m., approximately when core temperature reaches its maximum) (*32*) (see Supplement 2.2 for more detail).

## Results

A variety of measures (e.g., phase shift, variability, entropy) were examined for evidence of entrainment or any influence of the embedding rhythm on stability or attractor location (only entropy production was suggested for the main result; see Supplement 2.4 for the entropy calculation). With respect to experiment 1, behavioral performance (entropy) shows a maximum at 5:00 but has a more clearly defined minimum at about 17:00 in the day–night cycle, while the core body temperature rhythm shows a minimum at 5:00 but has a maximum at about 17:00 (Fig. 3.1). Regarding experiments 2 and 3, entropy was affected by the temporal locus during the circadian cycle, as well as by the introduction of a heated vest (Experiment 2) and an ice vest (Experiment 3); the effects of the thermal manipulation were not identical (see Supplement 2.2 for the temperature measure). Even if the same external temperature perturbations were given, the influence of the vest was negatively exaggerated (increasing entropy) at dawn, but positively exaggerated (decreasing entropy) in the evening (see Supplement 2.3 for more detail).
The estimated dynamics from the relative phase between two limbs, oscillatory coordination, was affected by the temporal locus during the circadian cycle (Fig. 3.2). Results at this biological scale correspond to a theoretical study which has shown that the rate of entropy production is changed when a new energy source is accessed via a nonequilibrium phase transition process (*33*).

## Discussion

The organism may convert the internal energy of itself so efficiently that it is able to produce what is physically possible (*34*). These results from Experiment 1, 2, and 3 reflect that accessing a new energy source differs as a function of the circadian cycle and that access can be manipulated by a temporary thermal manipulation (Fig. 4). Given the very widespread use of these features, what types of essential properties are involved in the extended emergent elementary dynamics between oscillators? According to researchers (*35*), the basic element of the coordination ($\phi$) is equal to $\phi$ ($x^{\theta_1 - \theta_2} = \phi$), and such an equation resembles the log base $x$ of $\phi$ is equal to $\theta_1 - \theta_2$. Then, with respect to the experimental results, the essential foundation of the symmetry dynamics between oscillators ($\phi$), the preferred elementary frequency of the individual segment of the individual segment of $x$ to the other relative phase from the intended phase, is nearly equal to the slightly asymmetric potential ($x^{\phi_{ave} - \phi_0} = \Delta\omega$). Thus, if this logic simply keeps going and the outcome is observed in terms of the approximate relative stability of this coordination dynamic, this logic will have "$x$ to the variation of the relative phase ($h\phi$) is equal to $\theta_1 - \theta_2 = \phi$ multiplied by $\Delta\omega$. That is identical to the [$x^{h\phi} = \phi \cdot (\Delta\omega)$], and this dynamic potential will be finally dependent on the [$\Delta\omega$ ($rad * x^{-1}$), $rad = radian$].

Pervasive interconnectedness—everything is connected with another thing or other things (*11*)—suggests that behavior is adapted to perceiving both the nested environmental properties and one's own nested behaviors—a union that organizes actions on surrounding circumstances (*36*). The observation of the direct and robust relationship between biological aspects (body temperature and motor synchrony) and an environmental process (circadian temperature cycle) may echo the adaptation of our system to the environment (*37*).

### Part 4: Approximated common property of the behavioral patterns

Under these observations, a system may be found to exhibit a variety of hitherto unobserved dynamical behavior, including cultural evolutionary characteristics and the coexistence of multiple search strategies. This study proposes that the features investigated are a particularly appropriate assumption in terms of obtain simplicity from complexity.

$$h(x) = m\big(s(x^{-1})\big), \tag{2}$$

$$x^{-1} \rightarrow s \xrightarrow{s(x^{-1})} m \xrightarrow{m\big(s(x^{-1})\big)}$$

The expression $h(x)$ represents our way of modeling that denotes wherever and whenever the evolutionary system is observed. This model takes $(x)$ and it inputs the $(x)$ into $(s)$ and gets out $\big(s(x)\big)$, and then the model inputs that into the $(m)$ and finally takes $m\big(s(x)\big)$. Going back to the insight by researchers (*10*), related to the minimal starting point for understanding any system as a function of interest at any level of interest (see Equation. 1), $(s)$, this includes the collection of all parts of elements $[C(s)]$ and comprises all other things influencing them $[E(s)]$, even the internal structure of the system $[S(s)]$. Thus, we input $(x)$ as arranged in an $(x^{-1})$ so that it somehow provides us with intuition about a system's fundamental properties $(s)$. This property and then inputs into function $(m)$ help us get to the point of understanding any system based on the system of interest at any level of interest.

As we take this composite function, it models a system that starts with individual segments $(x)$ as the input, and it shows the minimal starting point of the system $(s)$ that will be dependent on relative individual distance. Thus, how should a system that can be identified or predicted $(h)$ be related to how it depends on the individual segments in a given context. The fundamental properties demonstrated may be able to create a useful system dynamics reference, so that this functional pattern could be applied to various phenomena (Fig. 5).

### Closing remarks

Oddly enough, those curiosities could have the same overall answer: the creation of sophisticated functions from simple elements (*38*). Possible evidence for the association of this property is found when we compute an approximation of its sensitivity (*39*); the collective structure can show the possible entity as a function of the system's own unique set of behavior in a long-time limit $[\lambda = \lim_{n\to\infty} (r_k^n)^{1/n}]$. Here, logic can set the object's crucial variable $(r_k^n)$ that causes the different value in both dynamics and can be considered as simply $(1/n)$ arranged in certain rules (i.e., exponential). What this function reflects is the kinds of highly sensitive components that come from the rate of change which increased or decreased at a certain point. We observed this by measuring the contraction (stable system) or expansion (chaotic system) near the orbit of distance $[d(x_0, x_0 + \epsilon)]$, during the next iteration of distance $[d(f(x_0), f(x_0 + \epsilon))]$ (Fig. 5).

These are the simple rules regarding the primary characteristics of a dynamic system: they show that a simple function can serve as a basic principle that can be used to investigate various patterns (*40*). Regarding the assumption behind the simulation, if we reduce the system to individuals, there will be nothing left to study. Thus, an understanding of evolutionary behavior is not reliant on an individual lodged in a corner of the system. Rather, the basic conditions need to be in place and to be met; the complexities of the system will then fall into the same rules (*16*), and their patterns will be estimated (*1*).

**Acknowledgements:** This research was supported by an NSF Grant BCS-1344725 (INSPIRE Track 1), and the National Research Foundation of Korea (2016K2A9A1A02952017). All of them provided written informed consent to the study approved by the local ethics committee (SNUIRB No.1509/002-002) and conformed to the ethical standards of the 1964 Declaration of Helsinki (Collaborative Institutional Training Initiative Program, report ID 20481572).

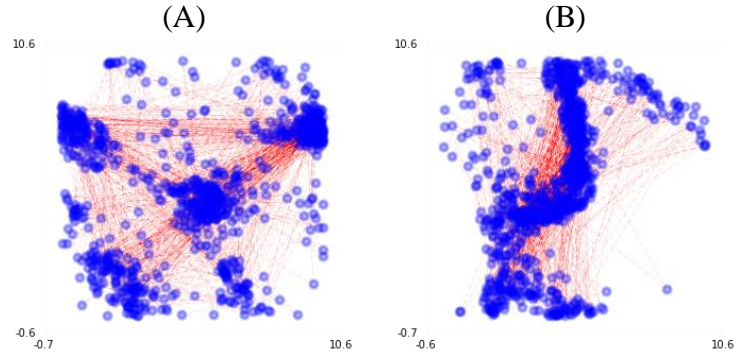

**Fig. 1.1.** Behavioral dynamics underlying social characteristics. Following the simulation, the left plot shows a displacement that separates individuals with a relative position structure controlled by the initial setting. This implies that although the pattern of individual behavior depends on a localized view of the initial conditions, a slight change in individual characteristics [IGT: individual's velocity up ($\vec{v}_i$) resulted in a loss of group heading ($\vec{v}_{avg}$)] underlying its social influence [ND: calculated from their social ties ($k = St$) multiplied by mutation rate ($u$)] has a remarkably diverging (A) or converging (B) effect on its displacement. The blue dots represent their position in an $x$, $y$ coordinate plane, and the red lines denote links (see Fig. S1 for more detail).

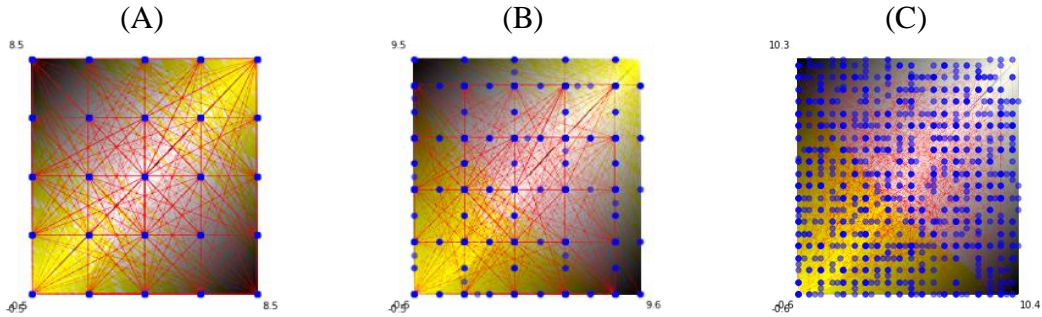

**Fig. 1.2.** Approximation of the evolution underlying interconnected interactions. The plots indicate that the patterns which occur correspond to the relative value. With certain defaults of their relativity, a slight change in the scalar value (social ties = St) derives dramatic impact at a certain point [A = St(0.55), B = St(0.56), C = St(0.57)]: blue dot = individuals, red line = links, background = density with symmetrical characteristics between individual and group heading. Notice that as the social ties increased (A ~ C), symmetrical charactersitics are biased to one side (C) (see Fig. S2 for more detail).

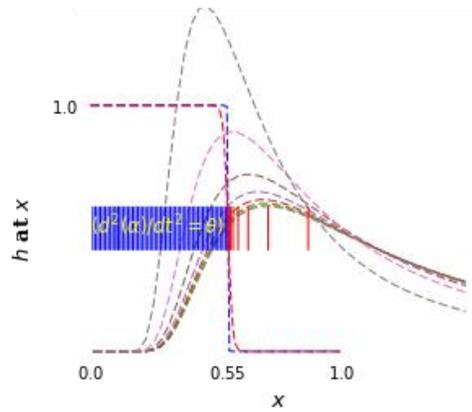

**Fig. 2.** Heuristics through the individual-based model. Based on the relativity defaults set by the model as an interconnected condition, the system becomes highly sensitive to small change in the scalar values (i.e., social tie) of individuals at a certain point. The horizontal axis denotes scalar (x = social ties in this simulation) from 0 to 1 and the vertical axis represents probability density at the scalar value (x). The plot suggests that if the individual fails to keep the trait (blue area = range from St 0.55 →) about the nearby individual, the displacement (red bars) will exponentially decay (dotted lines).

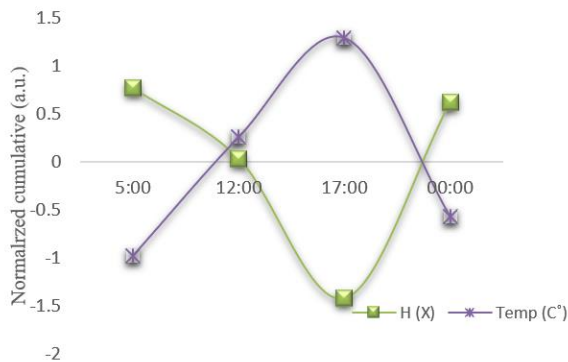

**Fig. 3.1.** Entropy production according to circadian cycles. Entropy features [H(x)] of the general tendencies in the normal condition [Temp(C°); see Supplement 2.2 for more detail on the entropy calculation]. Normalized = standard score (Z calculation), a.u. = arbitrary unit, H (x) = entropy production, Temp = temperature (Celsius), 5=5:00, 12=12:00, 17=17:00, 00=24:00. See Supplement 2.3 for further detail on the results.

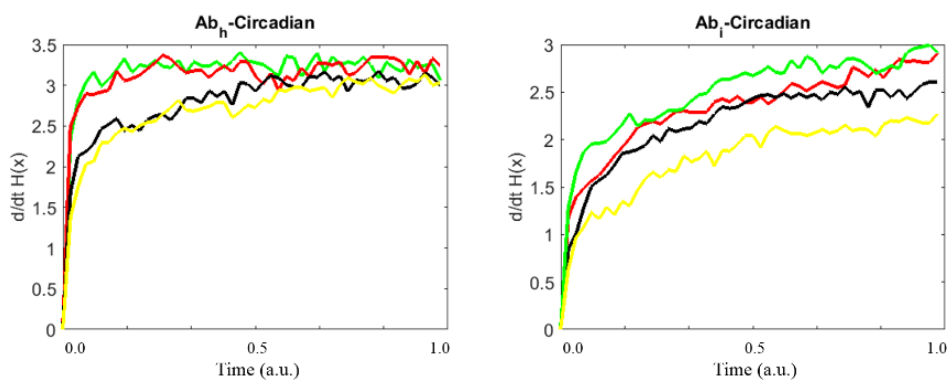

**Fig. 3.2.** Circadian and temperature perturbation dependent influences. Plots denote the estimated entropy forces according to the time series. The plot on the left denotes the entropy forces of the heat-based normal (Ab$_h$-circadian: red line = 5:00, black line = 17:00) vs. abnormal (green line = 5:00, yellow line = 17:00) conditions. The plot on the right denotes the entropy forces of the ice-based normal (Ab$_i$-circadian: red line = 5:00, black line = 17:00) vs. abnormal (green line = 5:00, yellow line = 17:00) conditions. See Supplement 2.3 for further detail on the results.

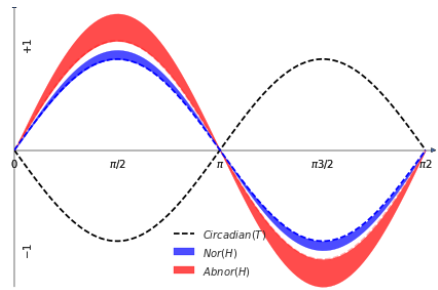

**Fig. 4.** Schematic illustration of the results according to the different experimental designs. The horizontal axis is the 24-hour circadian process as expressed by a sine function (pi/2 = 5:00, pi = 12:00, pi3/2 = 17:00, pi2 = 00:00), and the vertical axis is the optimized value of the state of the system with arbitrary units of -1 to 1. The blue line and shade (distribution) shows the observed normal states (H = entropy production) of the biological system according to the circadian temperature cycle. The red line and shade (distribution) denote the observed abnormal states (H = entropy production) in the perturbed circadian temperature conditions. Note that the investigation compared the biological values (relative phase) under an identical circadian condition (5:00 am, 5:00 pm), but different effects remain in terms of the stability of the system and the displacement (red area).

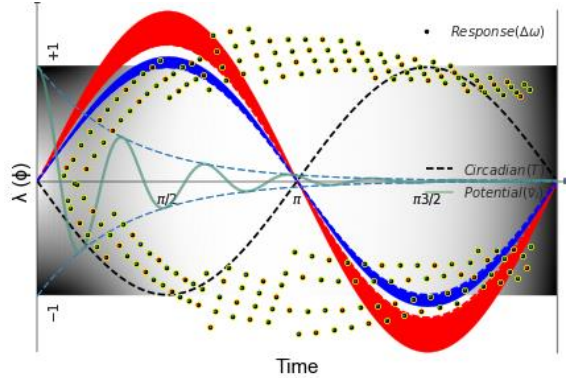

**Fig. 5.** Schematic illustration of the evolutionary understanding of the behavioral property. The plot represents the state of the system λ(φ) (one arbitrary cycle from -1.0 to 1.0) over time (horizontal axis). The green line indicates the damping force from the model 1 [decay a maintained value of $d^2(\vec{v})/dt^2$] between the focal individual and neighbors (or role individual) over time. The contour (black ~ white) represents the 24-hour circadian process as expressed by (pi/2 = 5:00 h, pi = 12:00, pi3/2 = 17:00, pi2 = 00:00) according to the optimized value of the system's state with arbitrary units of -1 to 1. The dotted lines show the observations from the model 2 (experimental results). The black line denotes the temperature (T) process according to the circadian cycle. The blue line and shade (distribution) shows the observed normal states of the biological system according to the circadian temperature cycle. The red line and shade (distribution) denote the observed abnormal states in the perturbed circadian temperature conditions. The dots surrounded by yellow colors (Response) denote plausible evidence for the association of this property. The crucial variable, which can intuitively be set in both dynamics, can simply be considered as the rate of change between the objects arranged according to a high-sensitivity rule [ $X_{n+1}$ is $r\chi_n(1 - X_n)$ ]. The results above describe a complex behavior with a divided phase (Δω) space in which areas of stability are surrounded by confusion. It implies that although their initial states are almost identical (in a comparison of the middle left area of the plot), the response becomes remarkably difference with iteration of $n$ times.

# Supplementary Materials for

**Evolutionary Understanding of the Conditions Leading to Estimation of Behavioral Properties through System Dynamics**

Chulwook Park*

Correspondence to: parkc@iiasa.ac.at

**This PDF file includes:**

Materials and Methods
Supplementary Text
Figs. S1 to S22
Tables S1 to S11
References (*1-30*)

# Contents

# 1. Model 1

The broader agenda of this supplement is to show the mathematical process behind the fundamental modeling mechanisms used. This information is based on spatially explicit mobility, where the individuals can move around their environment. The rules and processes in the artificially modeled structure describe an individual's homogeneous drives and also how to apply the mechanism.

### 1.1. Model detail

The agents are physically related to each other on some spatial representation allowing them to move anywhere in the space. The set of n-tuples of a real number, denoted by $\mathbb{R}^n$, is called n-spaces [$x = (x_1, x_2, \ldots, x_n) \in \mathbb{R}^n$]. A particular n-tuples in $\mathbb{R}^n$ is a point which called the coordinates, components, or elements of $x$. This is one of the standard ways in which the agents can continue to move in the space. The agents then move within the boundaries of the plane steering toward somewhere.

$$\vec{u} = \vec{a}_{avg} + \vec{b}_i, \qquad \vec{a}_{avg} = \|a_{avg}\| * \vec{d}_{avg}, \qquad \vec{b}_i = \|b_i\| * \vec{d}_i$$

where the $\vec{a}$ is the group's heading and $\vec{b}$ is each agent's ($i$) coherence toward the center of the group. The order ($\vec{u}$) is symmetric because all the agents are identical; thus individuals are naturally heading together in a certain direction, while at the same time maintaining a certain distance from each other as their inherited survival strategies (*1*).

The model, then, includes another operation with respect to the individual's current movement ($\vec{v}$) written simply by ($\vec{\omega} = \vec{u} + \vec{v}$). That new quantity of $\vec{\omega}$ is the sum of $\vec{u} + \vec{v}$, where the vectors stay away from the origin. The way we define this is that each vector represents a certain movement; a step with a certain distance and direction in space. If we take a step along the first vector of the $\vec{u}$, and then a step in the direction and distance described by the second vector of the $\vec{v}$, the overall effect is just the same as if we had moved along the sum of those two vectors to begin with.

$$\vec{\omega} = \vec{u} = \begin{bmatrix} 4 \\ 1 \end{bmatrix} + \vec{v} = \begin{bmatrix} 1 \\ 2 \end{bmatrix} = \begin{bmatrix} 4+1 \\ 1+2 \end{bmatrix} = \begin{bmatrix} 5 \\ 3 \end{bmatrix}$$

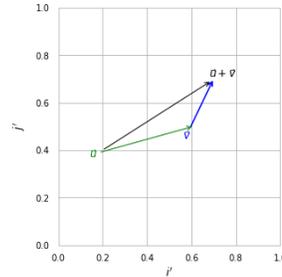

**Figure S1.** Schematic representation of the operation in $\mathbb{R}^2$.

The first quantity here (green line) has the coordinate ($\vec{u} = \begin{bmatrix} 4 \\ 1 \end{bmatrix}$), and the second quantity (blue line) has the coordinate ($\vec{v} = \begin{bmatrix} 1 \\ 2 \end{bmatrix}$). When we take the sum of the two quantities, we can see a four-step path from the origin to the tip of the second quantity: move 4 to the right and 1 up, then move 1 to the right again and 2 up. To rearrange these steps, first move 4+1 to the right, then move 1+2 up; the new quantity (black line) has coordinates 4+1 and 1+2 from the origin. Exhibited in this list-of-numbers conception is a matching up of their numerical terms, and an adding them both together ($\vec{\omega} = \begin{bmatrix} 4+1 \\ 1+2 \end{bmatrix} = \begin{bmatrix} 5 \\ 3 \end{bmatrix}$). With this fundamental process, especially, we see that the second quantity of $\vec{v}$ contains the individual's trait (agent-self interactions: agents can interact with themselves), based on its condition holds;

$$\vec{v} \to \vec{f} \to f(\vec{v}), \quad f(\vec{v}) = \begin{cases} \vec{v}\,(+)\; if\; f(\|v_i\|) < \|v_{avg}\| \\ \vec{v}\,(-)\; if\; f(\|v_i\|) > \|v_{avg}\| \end{cases}$$

where the function assumes that the attribute of the component is conditional upon the value yield in the other direction (-), if the length of the magnitude $\|v_i\|$ is greater than the other length ($\|v_{avg}\|$). This trait is implemented according to quantity as follows:

$$\vec{v}_i = \|v_i\| * \vec{d}_i, \qquad \|v_i\| = \sqrt{v_{i_x}^2 + v_{i_y}^2}, \qquad \vec{d}_i = \frac{\left( v_{i_x}, v_{i_y} \right)}{\sqrt{v_{i_x}^2 + v_{i_y}^2}}$$

$$\vec{v}_{avg} = \frac{1}{N} \sum_{i=1}^{N} \vec{v}_i = \|v_{avg}\| * \vec{d}_{avg}$$

where the $\vec{v}_i$ is the individual's velocity represented by the length of the individual's magnitude ($\|v_i\|$) with the direction of individual ($\vec{d}_i$). The $\vec{v}_{avg}$ is their average velocity which includes the entire population of N individuals' navigation. The result of the $\vec{v}_i$ and the $\vec{v}_{avg}$ produces a new quantity $\vec{v}_{inew}$ written simply as ($\vec{v} \to \vec{v}_{inew} = \vec{v}_i + \vec{v}_{avg}$). Going back to the conditional (IF) assumption with this individual quantity, if the object ($\vec{v}$) faces the parameters ($\vec{u}$) with the states of their quantities ($\|v_i\|$ and $\|v_{avg}\|$), the condition set produces an opposite direction (±) depending on its norm as follows:

$$f(\|v_i\|) < \|v_{avg}\| \quad \to \quad \vec{\omega} = \vec{u} + \vec{v}$$

$$f(\|v_i\|) > \|v_{avg}\| \quad \to \quad \vec{\omega} = \vec{u} + (-\vec{v})$$

where $\vec{\omega}$ is a new position vector of the individual updated by the inhereted trait $\vec{u}$ with the individual's current movement $\vec{v}$. This is a linear combination (or inverse transformation =180° counterclockwise), something that takes in inputs and spits out an output for each one;

$$\underbrace{\begin{bmatrix} x_{in} \\ y_{in} \end{bmatrix}}_{input} \to f \to \underbrace{\begin{bmatrix} x_{out} \\ y_{out} \end{bmatrix}}_{output} = \underbrace{\begin{bmatrix} 0 & 1 \\ 1 & 0 \end{bmatrix}}_{matrix} \underbrace{\begin{bmatrix} x \\ y \end{bmatrix}}_{vector} \to f \to \underbrace{\begin{bmatrix} 0 & 1 \\ 1 & 0 \end{bmatrix}}_{matrix} \underbrace{\begin{bmatrix} x \\ y \end{bmatrix}}_{vector} \; or \; \underbrace{\begin{bmatrix} 0 & -1 \\ -1 & 0 \end{bmatrix}}_{matrix} \underbrace{\begin{bmatrix} x \\ y \end{bmatrix}}_{vector}$$

$$\begin{bmatrix} x_{out} \\ y_{out} \end{bmatrix} = \begin{bmatrix} 1\vec{x} \\ 1\vec{y} \end{bmatrix} \; or \; \begin{bmatrix} -1\vec{x} \\ -1\vec{y} \end{bmatrix}$$

Imagine that every possible input vector multiplied by the matrix moves over to its corresponding output vector multiplied by the matrix (or inverse) without becoming curved and that the origin must remain fixed in place; what the coordinates are is determined by where each basis vector lands. For example:

$$\vec{\omega} = \vec{u} = \begin{bmatrix} 4 \\ 1 \end{bmatrix} + \vec{v} = \begin{bmatrix} 1 \\ 2 \end{bmatrix} = \begin{bmatrix} 4+1 \\ 1+2 \end{bmatrix} = \begin{bmatrix} 5 \\ 3 \end{bmatrix}$$

$$\vec{\omega} = \vec{u} = \begin{bmatrix} 4 \\ 1 \end{bmatrix} + \vec{v} = \begin{bmatrix} (-)1 \\ (-)2 \end{bmatrix} = \begin{bmatrix} 4+(-)1 \\ 1+(-)2 \end{bmatrix} = \begin{bmatrix} 3 \\ -1 \end{bmatrix}$$

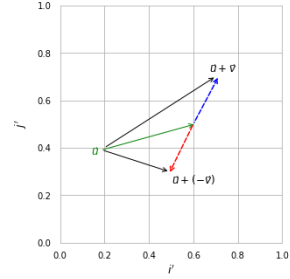

**Figure S2.** Schematic representation of the operation in $\mathbb{R}^2$.

This refers to the fact that the model's basic pattern of group behavior depends on the value of $\vec{v}$ with a localized view of the initial conditions of the randomly initialized point; increasing the individual's quantity $\vec{v}$ underlying the group's initial condition $\vec{u}$ causes their portrait to diverge (dotted blue line) or converge (dotted red line). Such a fundamental operation allows us to reach every possible point (not as an arrow but actually as a single point) in the plane, considering every possible linear combination that we can obtatin from the two dimensional quantities.

Based on these functions, the element $\vec{v}$ then contains a more detailed algorithm of how the individual's new position was implemented in which a subset holds that (i) the zero vector belongs to $\vec{v}$, and (ii) vectors and any muliplication of scalars is also in $\vec{v}$. Such a mechanism is dependent on quantity as follows.

First of all, the new position dynamic is augmented by the designated trade-off value (2) of the individual velocity-group heading [IGT: individual's velocity up ($\vec{v}_i$) resulted in a loss of group heading ($\vec{v}_{avg}$)] as follows;

$$\vec{v}_{iinew} = (1 - \|k\|) * \vec{v}_i + \|k\| * \vec{v}_{avg}, \qquad k \in [0,1]$$

where $\vec{v}_i$ is the velocity of each individual, $\vec{v}_{avg}$ is the average velocity about group heading, and the value of $k$ is a scalar that controls their trade-off. For example, the product of a $\vec{v}_i$ by a scalar $k$ is a vector $\|k\|\vec{v}_i$ with magnitude $\|k\|$ times the magnitude of $\|v_i\|$ and with direction $\vec{d}_i$, the same as or opposite to that of $\vec{v}_i$, according to whether $k$ is positive or negative [if $k = 0$ (null vector) $\|k\|\vec{v}$ has zero magnitude and no specific direction]. This means that every new position vector will be a combination multiplied by the scalar [$(c_1 = 1 - \|k\|), (c_2 = \|k\|)$];

$$c_1 \begin{bmatrix} 4 \\ 1 \end{bmatrix} \pm c_2 \begin{bmatrix} 1 \\ 2 \end{bmatrix} = \begin{bmatrix} c_1 4 \\ c_1 1 \end{bmatrix} \pm \begin{bmatrix} c_2 2 \\ c_2 2 \end{bmatrix} = \begin{bmatrix} x_1 \\ x_2 \end{bmatrix}$$
$$4 * c_1 \pm 1 * c_2 = x_1, \qquad 1 * c_1 \pm 2 * c_2 = x_2$$

and, the mechanisms obtain:

$$f(\|v_i\|) < \|v_{avg}\| \quad \rightarrow$$

$$\vdots$$

$$f(\|v_i\|) > \|v_{avg}\| \quad \rightarrow$$

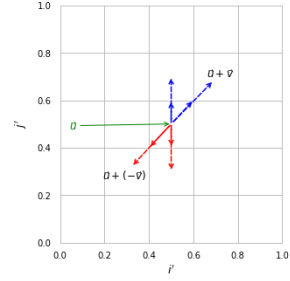

**Figure S3.** Schematic representation of the operation in $\mathbb{R}^2$.

where the different angles of the arrows are from $\vec{v}_{inew}$ and the spaces inside the arrows are from the multiplied scalar ($\|k\|, \ k \in [0,1]$). This gives us a global view so that we can conceptualize the list of quantities in a visual way and thereby simplify and clarify basic operational patterns.

The model then includes another characteristic for dealing with more or less distinct patterns of behavior. In social practice, just as individuals are more likely to change their decisions depending on the influences surrounding them (*3*), so the individual's new quantity ($\vec{v}$) holds its network characteristics (*4*) as given;

$$\vec{v}_{iiinew} = \left((1 - \|k\|) * \vec{v}_i + \|k\| * \vec{v}_{avg}\right) + \vec{v}_s, \qquad \vec{v}_s = \|v_s\| * \vec{d}_s$$

where $\vec{v}_s$ is a vector with a length $\|v_s\|$ and direction $\vec{d}_s$ as a function of the network density (ND). The network density arises from its social ties (*4*) based on the nodes (N = initial average exploration of the model) and calculated by its actual connection (AC) with the potential connection (PC) of the network.

$$\|v_s\| = \text{ND} = AC/PC, \qquad AC = (2 * t)/N, \qquad PC = N(N - 1)/2$$

where the network density ($\|v_s\|$) describes the potential connections in a network that are actual connections ($AC/PC$). The potential connection ($PC = N(N - 1)/2$) is a connection that could potentially exist between two individual regardless of whether or not it actually does. *This* individual could know *that* individual; *this* object could connect to *that* object. Whether or not they do connect is irrelevant when we are talking about a potential connection. By contrast, an actual connection ($AC = (2 * t)/N$) is one that actually exists. *This* individual does know *that* individual; *this* object is connected to *that* object. For example, in the living room of a house, the actual connections may represent one hundred percent of all the potential relationships. In contrast, on a public bus, they are likely to be quite low, relative to all the potential relationships, because the number of people who actually know each other on the bus (actual connection) must be low ($t$ =social ties).

Moreover, even on a public bus, any individual may connect to one another even if it knows nothing about the other individuals that it actually connects to (if the other offers the highest payoff). In the house, too, anyone can bring a guest into their living room. That may make the others modify the actual connection. These small linear contributions to their dynamics, and this structural instability can be interpreted as the network characteristics being influenced by the exploration rate ($k'$ = scalar), which corresponds to a mutation term in genetics given as:

$$\vec{v}_{iiiinew} = \left((1 - \|k\|) * \vec{v}_i + \|k\| * \vec{v}_{avg}\right) + \vec{v}_{ss}, \qquad \vec{v}_{ss} = \|v_{ss}\| * \vec{d}_{ss}$$
$$\|v_{ss}\| = [\|k'\|(1 - \|v_s\|) - 2\|k'\|\|v_s\|], \qquad k' \in [0,1]$$

where $k'$ controls how fast the transition function propagates in the network, and the new position vector considers its network density as another quantity $[(c_3 = \|k'\|)]$. For example, as previously:

$$c_1 \begin{bmatrix} 4 \\ 1 \end{bmatrix} \pm c_2 \begin{bmatrix} 1 \\ 2 \end{bmatrix} = \begin{bmatrix} c_1 4 \\ c_1 1 \end{bmatrix} \pm \begin{bmatrix} c_2 2 \\ c_2 2 \end{bmatrix} = \begin{bmatrix} x_1 \\ x_2 \end{bmatrix}$$
$$4 * c_1 \pm 1 * c_2 = x_1, \qquad 1 * c_1 \pm 2 * c_2 = x_2$$

The new position is then:

$$\begin{bmatrix} x_1 \\ x_2 \end{bmatrix} \pm c_3 \begin{bmatrix} 2 \\ 1 \end{bmatrix} = \begin{bmatrix} x_1 \\ x_2 \end{bmatrix} \pm \begin{bmatrix} c_3 2 \\ c_3 1 \end{bmatrix} = \begin{bmatrix} x_{11} \\ x_{22} \end{bmatrix}$$
$$x_1 \pm 2 * c_3 = x_{11}, \qquad x_2 \pm 1 * c_3 = x_{22}$$

In the presence of the network density, the system settles down into a state with a more pronounced increasing (or decreasing) mutation rate in every update step, and it yields:

$$f(\|v_i\|) < \|v_{avg}\| \quad \rightarrow$$
$$\vdots$$
$$f(\|v_i\|) > \|v_{avg}\| \quad \rightarrow$$

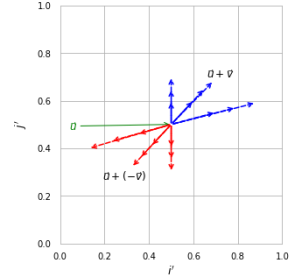

**Figure S4.** Schematic representation of the operation in $\mathbb{R}^2$.

With these implementations, instead of the widespread extension underlying the combination, the model proposes to adopt an existing possible interconnected relationship between the network and its movement characteristics. Let us think about simple interdependency between the two components (trade-off between individual velocity and group heading as an internal, network density multiplied by mutation as an external). If the individual's tendency is very remote from the group's purpose, its mutation in the system will not propagate to the individuals, or vice versa. For this application, the new position mechanism assumes that the social network characteristics ($\|v_{ss}\|$ = scalar) are a denominator applied by the index of difficulty ($id$ = scalar) as a numerator.

$$\vec{v}_{iiiinew} = [(1 - \|k\|) * \vec{v}_i + \|k\| * \vec{v}_{avg}] * (\|v_{id}\|/\|v_{ss}\|)$$
$$\|v_{id}\| = \frac{2D}{W}, \qquad \|v_{ss}\| = [\|k'\|(1 - \|v_s\|) - 2\|k'\|\|v_s\|]$$

where $\|v_{id}\|$ is the scalar as a function of the ratio between the two objects [2D = size of the trade-off ($k$) between two objects about $\vec{v}_i$ and $\vec{v}_{avg}$] divided by the width of the object [W = arbitrary value corresponding to the individual's size or reputation. For example, when $\|k\|$ is 0.1 applied to the $(1 - \|k\|) * \vec{v}_i + \|k\| * \vec{v}_{avg}$, the 2D becomes large (i.e., 0.8); on the other hand, when $\|k\|$ is 0.4, the 2D becomes small (i.e., 0.2)]. This leads to a simple interpretation linking the vector as follows;

$$c_4 \begin{bmatrix} x_1 \\ x_2 \end{bmatrix} = \begin{bmatrix} c_4 x_1 \\ c_4 x_2 \end{bmatrix} = \begin{bmatrix} x_{111} \\ x_{222} \end{bmatrix} = \vec{v}_{iiiinew}, \qquad \begin{bmatrix} x_1 \\ x_2 \end{bmatrix} = \vec{v}_{iinew}, \qquad c_4 = \|v_{id}\| / \|v_{ss}\|$$

and, the mechanisms yields;

$$f(\|v_i\|) < \|v_{avg}\| \quad \rightarrow$$

$$\vdots$$

$$f(\|v_i\|) > \|v_{avg}\| \quad \rightarrow$$

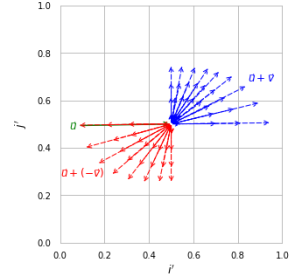

**Figure S5.** Schematic representation of the operation in $\mathbb{R}^2$.

Notice that the combination of these arrows refers to this model's fundamental feature. The characteristics of every rule and process that are applied (or will be applied) in this model must be within this functional dynamic. This provides an excellent way of conceptualizing many lists of individuals in a visual way, which can clarify patterns in mechanisms. It also shows a global view of what certain operations do to describe how an individual is being manipulated in space using numbers that can be run through a computation.

The model now considers an adoption probability which is given by an estimate of $\vec{v}_{inew}$ by the individuals. Indeed, as each individual may not know the exact value of the trait that has adopted the other's $\vec{v}_{inew}$, this model yields that they can estimate the value at every schedule of each generation via the comparison given.

$$p = [1 + e^{-\omega \Delta \pi}]^{-1}, \qquad \pi_r - \pi_f = \Delta \pi \big|_{\pi_r = role\ model}$$

where $p$ is the probability acceptance of the role model for imitation, $\pi_f$ is a payoff (velocity) of the focal individual, $\pi_r$ is a payoff of the role individual, $e$ denotes the exponential, and $\omega$ is the intensity of the selection ($\omega < 1$ = weak selection, $\omega \rightarrow \infty$ = strong selection). The focal individual imitates the strategy of the nearby role individual, comparing its new position vector (large $\Delta \pi$ = velocity difference large, small $\Delta \pi$ = velocity difference small), and then the focal individual chooses to imitate the strategy of the role individual.

The model applied this trait with three implementations ($\pi_1$, $\pi_2$, and $\pi_3$) with a different assessment of evolutionary patterns being expected from each of the above model mechanisms.

$$\pi_1 = \vec{u} \pm \left[ (1 - \|k\|) * \vec{v}_i + \|k\| * \vec{v}_{avg} \right]$$

$$\pi_2 = \vec{u} \pm \left[ \left( (1 - \|k\|) * \vec{v}_i + \|k\| * \vec{v}_{avg} \right) + \vec{v}_{ss} \right]$$

$$\pi_3 = \vec{u} \pm \left[ \left( (1 - \|k\|) * \vec{v}_i + \|k\| * \vec{v}_{avg} \right) * (\|v_{id}\| / \|v_{ss}\|) \right]$$

Given that the individuals naturally navigate together and maintain a certain distance as their inherited survival strategies, the equation contains a detailed algorithm of how the individual's new payoff (velocity) was implemented.

### 1.2. Table of the model variables

**Table S1. 1.** Variables of potential behaviors.

| | |
|---|---|
| Match | neighbors' average heading |
| Coherence | toward the center of the neighbors |
| Movement | individual displacement |

**Table S1. 2.** Model parameters.

| | | (Range of the value) |
|---|---|---|
| Parameters for defaults | | |
| Number of individuals | M | 1000 |
| Separation | D | $1 \sim 10$ |
| Vision (exploration) | E | $1 \sim 10$ |
| Velocity | V | $1 \sim 10$ |
| Parameters for internal (individual movement) characteristics | | |
| Individual velocity (from the velocity) | $v_i$ | $1 \sim 10$ |
| Group velocity (averaged the individual velocity) | $v_{avg}$ | $1 \sim 10$ |
| Individual-group velocity trade-off (systemic scalar) | $k$ | $0.1 \sim 0.9$ |
| Index of difficulty (from the velocity trade-off) | id | $0.1 \sim 0.9$ |
| Parameters for external (social network) characteristics | | |
| Node (from number of individuals) | n | $1 \sim 10$ |
| Social ties (systemic scalar) | t | $0.1 \sim 0.9$ |
| Mutation rate (systemic scalar) | $k'$ | $0.1 \sim 0.9$ |
| Selection intensity | $\omega$ | $1 \sim 10$ |

*1.3. Results supplements*

$$\vec{v}_{inew} = \vec{u} \pm \left[ \left( (1 - \|k\|) * \vec{v}_i + \|k\| * \vec{v}_{avg} \right) + \vec{v}_{ss} \right]$$

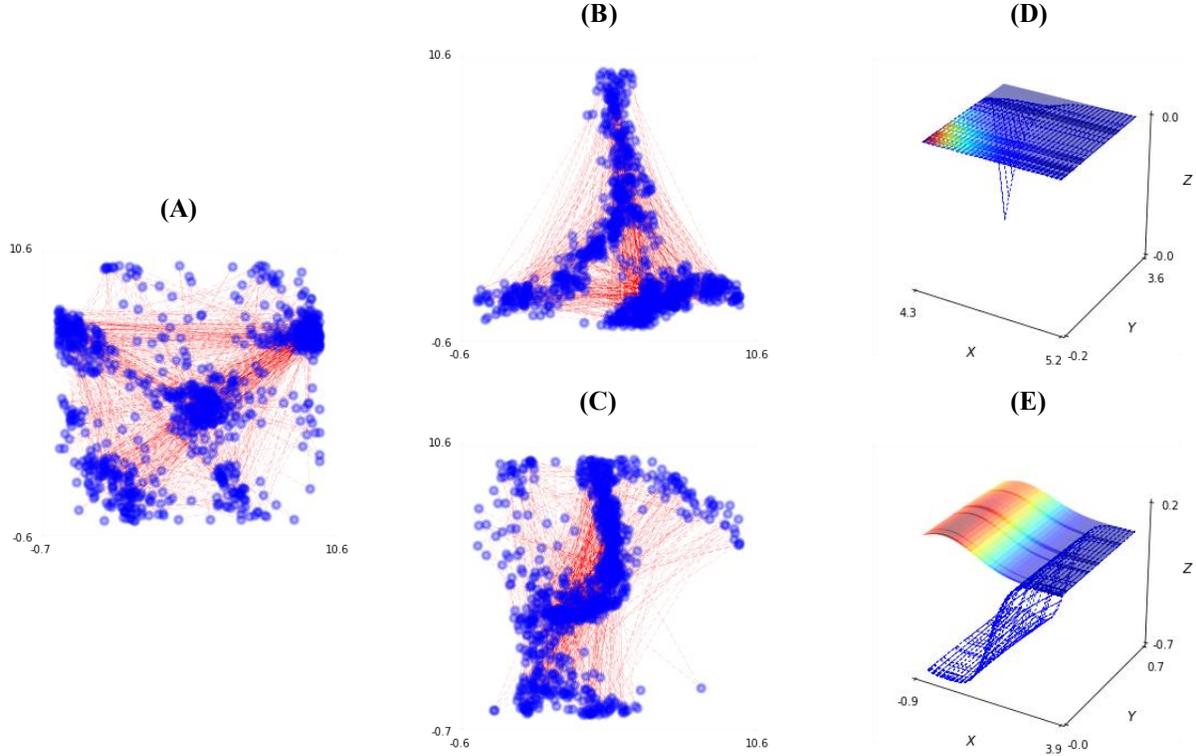

**Figure S6.** Behavioral dynamics underlying the model's simple movement rule (Equation = applied function: see Supplementary Materials 1.1 for more detail). From the simulation, the left plot (A) shows that a displacement separates the individuals, the relative position structure being controlled by the initial setting. It refers to the fact that although the pattern of the individual behavior depends on a localized view of the initial conditions, a slight change in the individual movement characteristics underlying its individual-group trade-off has a remarkably diverging (or converging) displacement. Blue dots represent their position in an x, y coordinate plane, and the red lines denote links. The plots of (B) and (C) show the social influence based on the network characteristics: the plot of the upper (B) was from the social ties multiplied by the mutation rate low with its density function in $\mathbb{R}^3$ (right side); plot of the bottom (C) was from social ties multiplied by the mutation rate high with its density function in $\mathbb{R}^3$ (right side).

$$\vec{v}_{inew} = \vec{u} \pm \left[ \left( (1 - \|k\|) * \vec{v}_i + \|k\| * \vec{v}_{avg} \right) * (\|v_{id}\| / \|v_{ss}\|) \right]$$

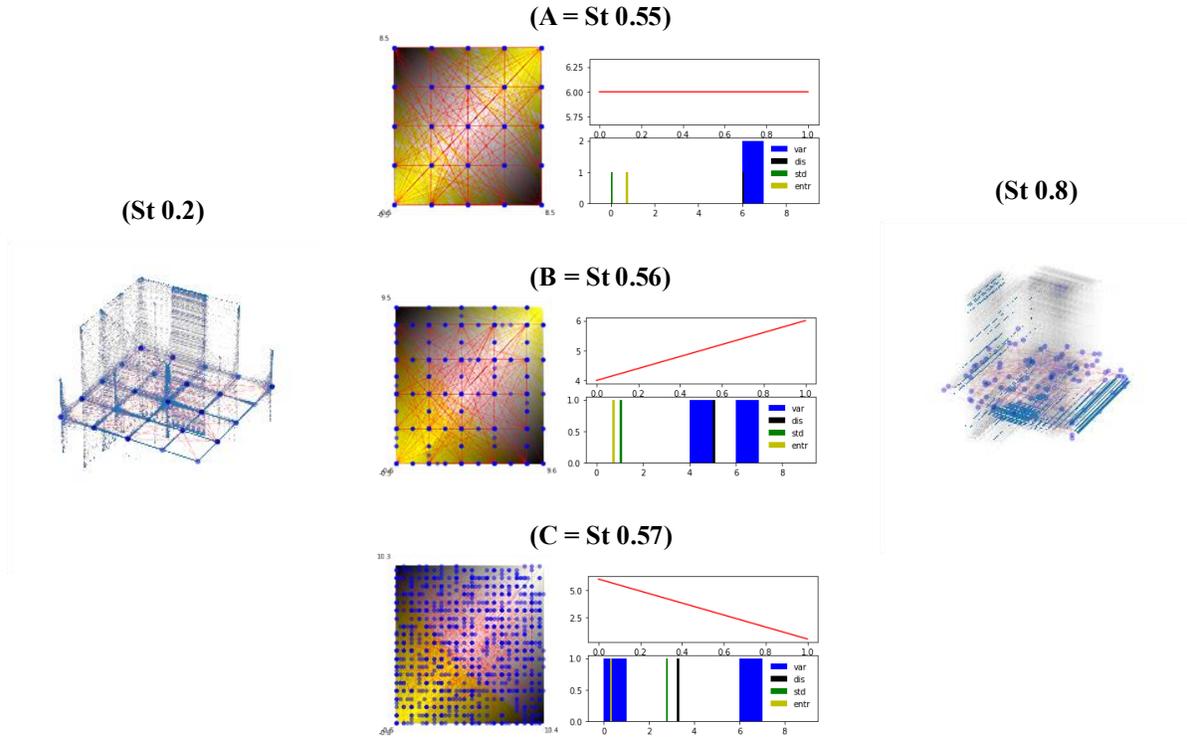

**Figure S7.** Assessment of the strategy and evolution underlying its interconnected interactions (Equation = applied function: see Supplementary Materials 1.1 for more detail). The plots indicate that the changes occur at a certain point.

**In the middle plots** (A, B, C), blue dots = individuals; red line = links; background = density with symmetrical characteristics between individual and group heading; black bar denotes average displacement (red line) of the individual's position; blue bars indicate its minimum and maximum clustering variance; green bar = standard deviation; and yellow bar = entropy. Note: the system becomes highly sensitive to tiny changes in the individuals' social ties (considered at St. 0.55 in this simulation) at which the model is set as an initial value.

**In the figure left**: Social ties = 0.2, mutation rate = 0.5; IGT = 0.1, ID = 0.8; blue dot = individuals; red line = links; background = density in $\mathbb{R}^3$

**In the figure right**: Social ties = 0.8; mutation rate = 0.5; IGT = 0.1, ID = 0.8; blue dot = individuals; red line = links; background = density in $\mathbb{R}^3$

## 2. Model 2

The supplementary information reported in this section is analytical information showing how interacting cyclic processes account for the emergence of new entities. It investigates whether unintentional coordination provides an environmental rhythm within an individual's field of view, and will explain whether the dynamics of bimanual coordination is influenced by an overarching temporal structure that is irrelevant to the task.

### 2.1. Model detail

***Circadian rhythm of temperature (external source):*** The core cycles of a biological system are influenced by temperature, with 24-hour light-dark oscillation (called circadian rhythm), as well as by biochemical, physiological, or behavioral processes that persist under constant conditions with a period length of ~24 hours (*5*). Presumably, due to inputs to the thermoregulatory centers from the body core (*6, 7*), the circadian rhythm of biology shows a minimum at 5:00 (when core body temperature is rising most rapidly) but has a more clearly defined maximum at about 17:00 in the daylight (when core body temperature is falling most rapidly) cycle (*8, 9*).

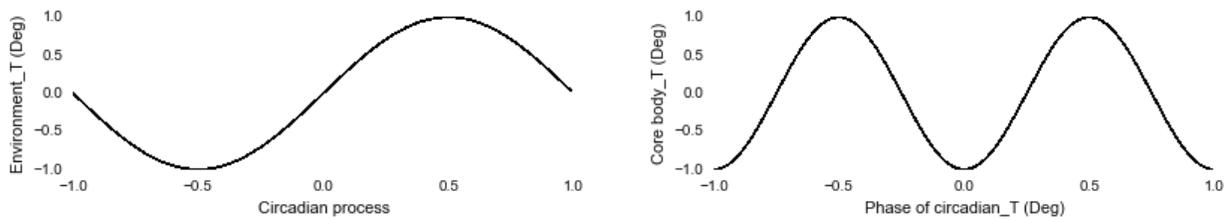

**Figure S8.** Representation of the circadian rhythm. Left = circadian process oscillation; right = temperature process oscillation between the circadian temperature (horizontal axis) and the body temperature (vertical axis). Note: this is a normalized rhythm despite the fact not all rhythms are identical. Our core body temperature is roughly linked to this cycle, with various hormones being released at certain stages during the rhythm because our body temperature reflects energy levels.

This circadian change (in core temperature) is most likely due to the rhythmic input from the suprachiasmatic nuclei (SCN) acting upon the hypothalamic thermoregulatory centers and altering the thresholds of cutaneous vasodilatation and sweating (*10*). Specifically, melatonin appears to contribute to this change, as its rate of secretion increase in the evening, and this increase promotes a fall in body temperature via cutaneous vasodilatation (*11*).

As the information is accessible, people are familiar with how such a process can fluctuate and how it can be explained by the interaction between the internal (homeostatic) and the external (circadian) situations (*12*). There is ample evidence of the effect of the ecological climate on various aspects of the process (*13*). Heat exchanged with the environment by means of convection and radiation allows a gradient to be formed between the body core and temperature (*14*). The rhythm in the core temperature produced by this change is generally promoted by other rhythms, including the body clock, sleep, and physical and mental activity, raising the possibility that the disruption of circadian rhythms can contribute to complications in the human system (*15*). These changes in the interior temperature in the body, as opposed to the peripheral (core

temperature)—both in animals and humans—are mainly due to circadian rhythmic changes in the rates of ecological impacts (*16*).

However, given that precise control of the internal substance (SCN) as a generator of biological circadian rhythm is unclear, the circadian rhythm of the core body temperature appears to be generated by periodic variation in heat production and heat loss (*17*). For instance, changes in heat loss via convection and radiation are primarily caused by variations in skin blood flow, with consequent changes in skin temperature. In particular, when the subject is performing mild activities, where a decreased temperature is not matched by a thermal load, it has been shown to be very effective in describing the thermal responses to activity carried out at different times of the day (*16*).

When one considers the submaximal activity changes following a brief period, say, at a certain temperature level (*18*), one may see initially that the response to the same amount of moderate activity in the minimum circadian rhythm differed from that in the maximum circadian rhythm. The mechanisms responsible for these different temperatures of the core and musculature during daylight cycles, as a result of normal or non-normal ambient temperatures, will alter a range of performance factors, including the thermoregulatory response to activity. These results fully substantiate the predictions based on the hypothesis describing a circadian rhythm in thermoregulatory responses and indicate that this hypothesis applies to biological adaptation regarding certain ecological variables.

***Elementary coordination of the HKB model (internal source):*** Formation and retention refer to propriospecific information about the states of the muscular-articular links, and the dynamical criteria of the stability pattern constrain the patterns or characteristics. To be specific, let us consider a qualitative physical system such as stiffness, damping, and position over time in a dynamical mass-spring system as given.

$$f(t) = mx'' + bx' + kx \qquad (1\text{-}1)$$

Here, $m$ is mass, $b$ is friction, and $k$ denotes the stiffness. The variable t is time, $\chi$ denotes the position, $\chi'$ is velocity, and $\chi''$ represents acceleration. In physics, because damping is produced by a process that dissipates the energy stored in the oscillations, the interplay between input and damping approaches a stationary fixed point in the long-time limit.

$$mx'' + bx' + kx = 0 \qquad (1\text{-}2)$$

Such systems possess a static equilibrium point, which is called a point attractor (*19*).

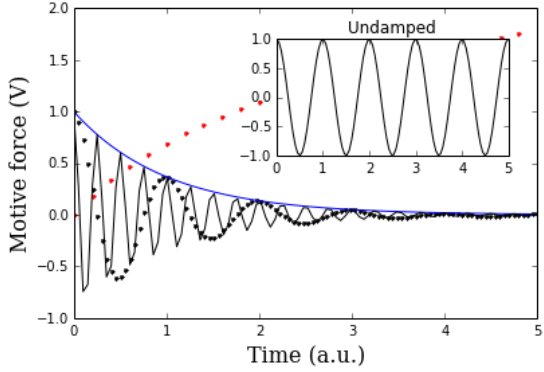

**Figure S9.** Simulation of the different mass-spring attractors. The damped exponential decay of the dotted equals $\cos(2\pi t)\exp(-t)$, and for the solid line $\cos(8\pi t)\exp(-t)$. Log lines indicate an embedded invariant property in terms of the relation between systems' attractor (solid blue line) and damping potential (dotted red line). V = volts, a.u. = arbitrary unit. Note: inserted plot denotes an undamped case.

The property of this dynamic has been applied not only to a physical system but also to descriptions of the human neuromuscular level (*20*). This function involves an investigation of the intact movement of a limb oscillator in terms of muscle-joint kinematic variations (kinematic position, velocity, acceleration) over time. When we are asked to swing two limbs comfortably, this can be characterized by the pendulum's dimension (*21, 22*), namely, simplifying the point attractor while restricting it to certain domains of phase space.

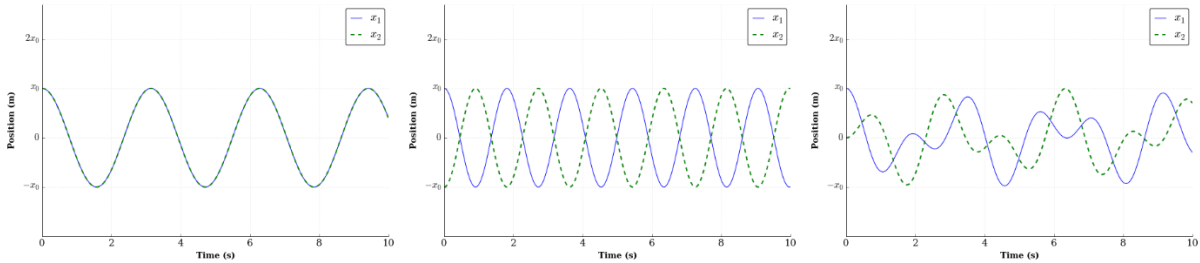

**Figure S10.** Synchronous diagrams of the possible point attractors: left = synchronized almost in-phase, with a phase difference x2 − x1 ≈ 0 and, in the anti-phase condition = middle, when x2 − x1 ≈ π, or it does not match the x2_initial input 0.5, x1_initial input 0.0 = right.

$$\theta 2 - \theta 1 \approx 0 \tag{2-1}$$

$$\theta 2 - \theta 1 \approx \pi \tag{2-1}$$

In this equation, with the phase difference, $\theta 2 - \theta 1 \approx 0$ denotes a condition of nearly synchronized in-phase, and $\theta 2 - \theta 1 \approx \pi$ indicates that this in an anti-phase. The observed relative phase or phase relation (ϕ) between two oscillators at ϕ ≈ 0 deg (in-phase), or ϕ ≈ 180 deg (anti-phase) have been modeled as the point attractors in our limb system, as they are purely stable patterns (*23*).

In the observed relative rhythmic segments patterns, the in-phase ϕ = 0 condition is more stable than the anti-phase ϕ = π condition. Inspired by a number of studies on the 1:1 frequency locking of the left- and right-hand phase defined as ϕ = $(\theta_L - \theta_R)$—the difference between the left (L) and right (R) phase angles (ϕ)—has led to the identification of important invariant human system features (*24*).

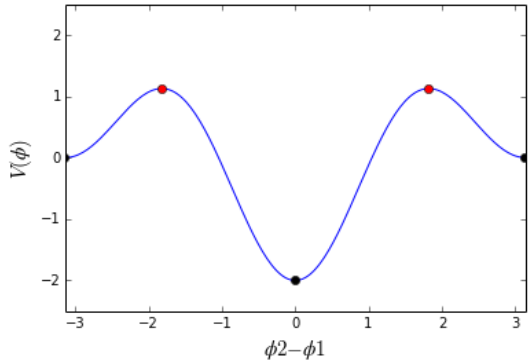

**Figure S11.** Reflection of the prototype (25). The vertical axis denotes the energy of the function at each averaged relative phase. The horizontal axis indicates the averaged relative phase between two limbs from in-phase 0 to anti-phase 180 (-180). At the local point of 0 and 180 (-180), the function is close to those minima (attractors = black balls) and at local point around 90 (-90), the state is close to the maxima (repellors = red balls).

$$V(ϕ) = −α\,cos(ϕ) − b\,cos(2ϕ) \qquad (3)$$

In this equation, ϕ is the phase angle of the individual oscillator. In addition, $α$ and b are coefficients that denote the strength of the coupling between the two oscillators. The relevant regions of the parameter space allow the potential V(ϕ); the negative signs in front of the coefficients simplify the equation of motion. A relative 1:1 frequency-locked coordination phase [V(ϕ)] is determined by the differences between the continuous phase angle [$−α\,cos(ϕ) −$ b cos(2ϕ)] of the oscillator's two components: the stability of the point attractor can be varied by varying the pendulum's dimensions (*21*).

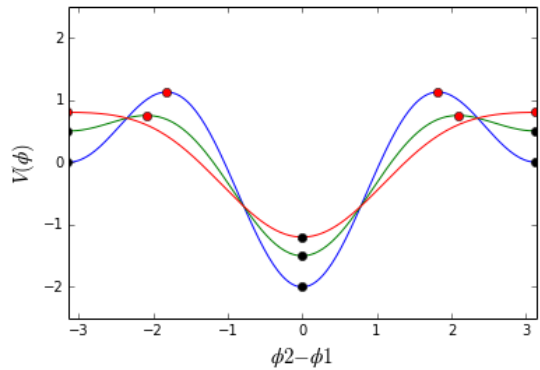

**Figure S12.** Reflection of a potential function. The blue line = the vertical axis, which denotes the energy of the function at each averaged relative phase. The horizontal axis indicates the averaged relative phase between two limbs from in-phase 0 to anti-phase 3.14 = 180 (-3.14 = -180). At the local point 0 and 180 (-180), the function close to those minima (attractors, black balls) and at local point around 90 (-90), the state is close to maxima (repellors, red balls). The red and green lines denote the variation of the potential functions. Depending on the two components' preferred frequencies [$−α\,cos(ϕ) −$ b cos(2ϕ)], the intrinsic dynamics of the potential function of V is determined. Black balls symbolize stable states (attractors) and red balls correspond to unstable states (repellors). While all states of the systems become less stable due to a higher energy level,

the balls at the relative phase point of 0 (in-phase) remain; these balls in the anti-phase condition disappear because the system is shallower.

This function indicates that the minima of the potential are located at φ=0, and that $\phi = \pm\pi$ (25). Given this scenario, the function can be estimated in terms of how the potential will change in shape, as the control parameter (energy cost) increases. Based on the observed mechanism for the point attractor with a simple function, the present study proposes the in-phase bimanual rhythmic coordination synchrony pattern as a particularly well-suited physical model. This allows a useful reference for system stability coordination tasks in which this functional pattern can be applied to all human movement, muscles, and even a neural network. Actual intersegmental coordination, however, is additionally shaped by the contingencies of adjusting to environmental vagaries. How these extrospecific requirements and information types are incorporated into the physical stability patterns can be assumed by the level of symmetry coordination (26). In order to harmonize the effects of motor stability toward environmental symmetry, this study investigates the following elaboration.

***Symmetry breaking in bimanual coordination dynamics:*** The potential [V(φ)] extends the described assumption in terms of the difference between the uncoupled frequencies of bimanual rhythmic components:

$$\Delta\omega = (\omega_L - \omega_R) \tag{4}$$

where $\omega$ is the preferred movement frequency of the left ($\omega_L$), right ($\omega_R$) individual. If the relative phase between $\omega_L$ and $\omega_R$ were equal ($\Delta\omega = 0$), this pattern would be assumed to be a perfectly identical symmetry. However, the preferred movement frequencies of the individual oscillators in in-phase are large (i.e., function: b/a=0.5, detuning=-0.5, or detuning=-1.5), the expected stability of the rhythmical limb oscillation dynamics become greater than equal.

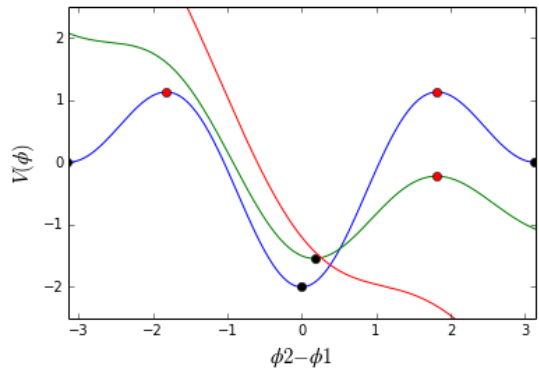

**Figure S13.** Preferred movement frequencies of the individual oscillators: blue line denotes the same symmetry; green line denotes the large different symmetry (function: b/a=0.5, detuning=-0.5); redline denotes even larger different symmetry (function: b/a=0.5, detuning=-1.5).

Such phenomena of the symmetry breaking must be another fundamental feature of the coordinative system (27). From this dynamic, a different noise of the underlying subsystems (neural, muscular, and vascular) can be estimated around an equilibrium point, and this might conceptualize the model when it comes to making operational definitions of each category in which the model has to consider the variability of the relative phase frequencies between two limbs:

$$\dot{\phi} = \Delta\omega - \alpha\,cos(\phi) - b\,cos(2\phi) + \sqrt{\varrho\xi_t} \qquad (5)$$

The estimation of two oscillators' relative phase ($\dot{\phi}$) is captured by the parameter ($\Delta\omega$) of the preferred movement frequency of the individual segment [$\alpha\,cos(\phi) - b\,cos(2\phi)$] with the noise ($\sqrt{\varrho\xi_t}$). Given the equation of the preceding model (grouped as the kinematics of motor stability according to the coordination task of synchronization), such a term has been used to capture purely functional dynamics regarding the equilibria and is confirmed usually as in the time and temporal difference between an oscillating limb (28).

Researchers (28), conducting experiments in handedness, advanced the elementary coordination dynamics of Eq. 5. They added two add (sine) terms for the coefficients, whose signs and magnitudes determine the degree and direction of asymmetry, as follows;

$$\dot{\phi} = \Delta\omega - [\alpha\,sin(\phi) + 2b\,sin(2\phi)] - [c\,sin(\phi) + 2d\,sin(2\phi)] \\ + \sqrt{\varrho\xi_t} \qquad (6)$$

Here, $\dot{\phi}$ indicates a coordination change. $\Delta\omega$ refers to a symmetry breaking through frequency competition between two limbs. [$\alpha\,sin(\phi) + 2b\,sin(2\phi)$] denotes a symmetric coupling defined by relative phase of 0 and $\pi$ attractors (this form of the term could be derived as the negative gradient potential V with respect to $\phi$); and the [$c\,sin(\phi) + 2d\,sin(2\phi)$] terms means added asymmetric coupling attractors with the stochastic noise $\sqrt{\varrho\xi_t}$. This extended equation refers to the fact that the emergent elementary dynamics between limbs or limb segments was governed by a slightly asymmetric potential of the [$c\,sin(\phi) + 2d\,sin(2\phi)$]. That suggests extended collective dynamics of the inter-segmental rhythmic coordination of the periodic components.

***Thermoregulatory symmetry breaking of the elementary coordination:*** Inspired by the complementary symmetric and asymmetric influences, the described model was applied to investigate the difference between the coupled or uncoupled frequencies of the temperature-rhythmic components between the core body and circadian cycles.

$$c = circadian \text{ temperature cycle}$$
$$d = core\ body \text{ temperature cycle}$$

where $d$ is the preferred rhythmic frequency of one (the homeostasis cycle) and another ($c$ = circadian cycle) individual. Whereas b/a determines the relative strengths of the fundamental in-phase equilibria, small values of $c$ and $d$ break the symmetry of the elementary coordination dynamics while leaving their essential coupling characteristics.

$$|c\ and\ d| > 0$$
$$|c\ and\ d| \approx 0$$

In this proposed assumption, the coefficient of the $d$ should be more important, producing the empirically observed perturbation in the equilibrium phase state, and then the $c$ should be set to zero without loss of generality, given that we cannot manipulate the environmental circadian cycle. As one can see, if the coupling between $d$ and $c$ is strong ($|c\ and\ d| \approx 0$), this pattern would potentially be expected to be in perfectly corresponding symmetry with the environmental requirement. However, the preferred rhythmic coupling of individual oscillators in an in-phase condition becomes a difference ($|c\ and\ d| > 0$), and thus the expected stability or variability of the rhythmical-component oscillation dynamics will become greater than equal. Given the

preceding assumption (grouped as the kinematics of motor stability according to the coordination task of synchronization), the equation was extended to a novel task in which there are different sources of symmetry breaking through thermal variables, as information has not yet been made available about the effects of bimanual dynamics in instruction on circadian temperatures.

$$\dot{\phi} = \Delta\omega - [\alpha\,sin(\phi) + 2b\,sin(2\phi)] - [c\,sin(\phi^{°C}) + 2d\,sin(2\phi^{°C})] \\ + \sqrt{\varrho\xi_t} \tag{7}$$

In this equation, in the bimanual 1:1 rhythmic coordination performed at different coupled frequencies, the symmetric coupling coefficients will be not the same. There will be an increase in detuning ($\Delta\omega$) and a decrease in the relative strengths of the attractors at 0 and $\pi$. However, when it comes to our limiting case of $\Delta\omega = 0$ on the approximately identical symmetry temperature parameters (core body and circadian cycle), what should we expect? The final estimation between the relative phases of two oscillators ($\dot{\phi}$) will be captured mainly by the parameter of the asymmetric thermoregulatory coupling $[c\,sin(\phi^{°C}) + 2d\,sin(2\phi^{°C})]$ with noise ($\sqrt{\varrho\xi_t}$).

From this dynamic, the different noise types of the underlying subsystems (neural, muscular, and vascular) around an equilibrium point were able to be estimated, suggesting that such phenomena related to symmetry breaking may be yet another remarkable feature of the coordinative system.

In sum, this experiment was required to have a condition of in-phase ($\phi = 0$) oscillated simultaneously at the 1;1 frequency locking (same tempo). The same goal using the functional symmetry dynamics of different effectors will be influenced by the asymmetric thermal regulation symmetry breaking through both circadian temperature cycles. Namely, the effect of one of the contralateral homologous relative limbs phase might be not identical to the impact of the others. The expected stability pattern, from intuition given a different motor, appears to allow the biological symmetry dynamic to be understood in the ecological context. This attunement to the circadian temperature approach implies an emergent property of the system.

### 2.2. Experimental design

**Design for experiment 1:** The present experiment was designed to verify whether the ecology influences physical systems. In experiment 1, to obtain the rate of motor synchrony depending on environmental cycles, the data were subjected to an analysis of variance comparing normal day-night temperature effects (four levels of circadian rhythm: 12:00 am, 5:00 am, 12:00 pm, 5:00 pm).

**Table S2.** Data collection for experiment 1: eight participants, four circadian points, six trials at each circadian point.

| Condition | Participants (N) | Circadian points | Trials at each circadian | Task/rest (min) |
|---|---|---|---|---|
| Normal | 8 | 5:00 12:00 17:00 00:00 | 6 | 1/5 |

Note: Participants were asked to swing their limbs in-phase at different anatomy points [192 datasets (three levels: wrist, elbow, and shoulder)], but only wrist joint data (64 set) were used for analysis (see Supplement 2.6 for more detail). The duration of each trial was 1 minute and there was a 5-minute rest interval between trials.

It refers to the response to a question "Does an ecological feature influence our system?" In-phase bi-manual coordination synchrony serves as a dependent variable, according to independent variables of the normal circadian temperature cycles.

***Design for experiment 2 and 3:*** In experiments 2 and 3, with respect to the question, "How does our system adapt to a regular or irregular thermal structures?" normal and abnormal day–night circadian temperature effects were compared at dawn (5 a.m.) and dusk (5 p.m.), given that our core temperature reaches its maximum at around 5 p.m. and its minimum at approximately 5 a.m.

**Table S3.** Data collection for experiment 2: two conditions, 8 participants, two circadian points, six trials at each circadian point.

| Condition | Participants (N) | Circadian points | Trials at each circadian | Task/rest (min) |
|---|---|---|---|---|
| Normal | 8 | 5:00 17:00 | 6 | 1/5 |
| Abnormal (heat based) | 8 | 5:00 17:00 | 6 | 1/5 |

**Table S4.** Data collection for experiment 3: two conditions, eight participants, two circadian points, six trials at each circadian point.

| Condition | Participants (N) | Circadian points | Trials at each circadian | Task/rest (min) |
|---|---|---|---|---|
| Normal | 8 | 5:00 17:00 | 6 | 1/5 |
| Abnormal (ice based) | 8 | 5:00 17:00 | 6 | 1/5 |

Note: Participants were asked to swing their limbs in-phase at different anatomy points [192 data set (three levels: wrist, elbow, and shoulder)], but only wrist joint data (64 set) were used for analysis. The duration of each trial was 1 minute with a 5-minute rest interval between trials.

Hence, in-phase bi-manual coordination synchrony serves as a dependent variable for two independent variables, as follows: two levels of circadian rhythm * two levels of thermal variable manipulation (experiment 2 = increasing, experiment 3 = decreasing).

***Apparatus and procedure:*** In-phase coordination without detuning was performed while each subject was seated in a chair holding a pendulum vertically without occluding their vision. The pendulums used here were two standard wooden rods (85 g, 1 m in length, 1.2 cm in diameter) with DC potentiometers attached. A 200 g weight was positioned 30 cm from the bottom of the rods. Each participant was asked to grasp the pendulum firmly 60 cm from the bottom so that the pendulum would not slip out of their hands, and not to rotate their finger joints. Their forearms were fixed voluntarily so that the pendulum motion was restricted to the sagittal parallel plane and the joint vertical axes (i.e., each oscillation pertained to only one joint, with the other joints being held immobile).

For experiment 1: the sessions were tapped into the ongoing circadian rhythm, focusing on its thermal structure. This involved four temperature (normal) embedding cycles (12:00 am, 5:00 am, 5:00 pm, 12:00 pm). Participants had to present four times with six trials each time (one participant, 1 minute, 24 trials = 6 [trials] *4 [circadian points]). Each trial block lasted for 1 minute with a rest of 5 minutes. The participants had received instructions about the preferred pendulum movements to establish in-phase 1:1 frequency locking at a 1.21 s metronome beat (this period was chosen because it corresponded to the natural period of the pendulum system) without concern over amplitude or frequency (*29*). A small amount of experience was provided to help avoid problems in complying with the session requirements, though instruction before the experiment. During the actual trial, no feedback was given, and the participants were not allowed to report except when a problem arose. If the participants accidentally moved a joint that was supposed to be voluntarily fixed, the data from that trial were not analyzed, and the trial was repeated at a later time. The three different oscillation joints were used in random order.

For experiment 2 and 3: the sessions were introduced for short-term thermal variable manipulation involving two conditions (normal and abnormal), two temperature-embedding cycles (5:00 am and 5:00 pm—the lowest / highest peak of the circadian rhythm of core temperature with skin capacitance). Each trial block lasted for 1 minute with a rest of 5 minutes. This trial and rest time are related to the maintenance of the thermal capacity of the body. In the natural session (normal temperature), participants had received instructions about the preferred pendulum movements, as designed for the previous procedure. For the perturbed condition (abnormal temperature: artificially decreasing/increasing temperature), participants donned a heat (or ice) vest for 30 minutes that would increase/decrease their core temperature. Their core temperature was checked using a monitoring system (3M$^{TM}$ Bair Hugger$^{TM}$ Temperature Monitoring System Control Unite and Sensor). After 30 minutes, after verifying the participants' abrupt temperature change, the data were collected for each trial in the same manner as in the natural condition setup. As there were concerns over the body temperature change, the participants' temperature was checked during each of the rest periods. The intermittent movement was chosen in a short time (within 30 minutes) taking into consideration the participant's body temperature capacities which were modified by pre-applied (heating or cooling) exogenous temperature (*10*). The entire session for each block did not last more than 30 minutes.

***Temperature measure and analysis:*** There is a state of unit $x$ at time t [$x(t)$] which denotes the ensemble average of the system. Collective behavior emerges from a homogeneous state when a parameter makes a transition from $x = 0$ to $x \neq 0$. In order to reflect this transition in terms of an external parameter, the temperature is such that all the terms of the parameter that can be considered as perturbations are considered in this case.

$$\epsilon\ (\chi, t) = \mathrm{T}\ (\chi, t) -\ \mathrm{T}_0(t) \qquad (16)$$

Here, $\epsilon$ is considered as the control parameter, while is T fixed. However, it is useful to capture the effective control parameter $\mathrm{T}_0$ to measure the distance of $\epsilon$ with respect to T. That is, the $\mathrm{T}_0 < \epsilon$ state is unstable (ice-vested and heat-vested), whereas the $\mathrm{T}_0 > \epsilon$ state is stable (normal).

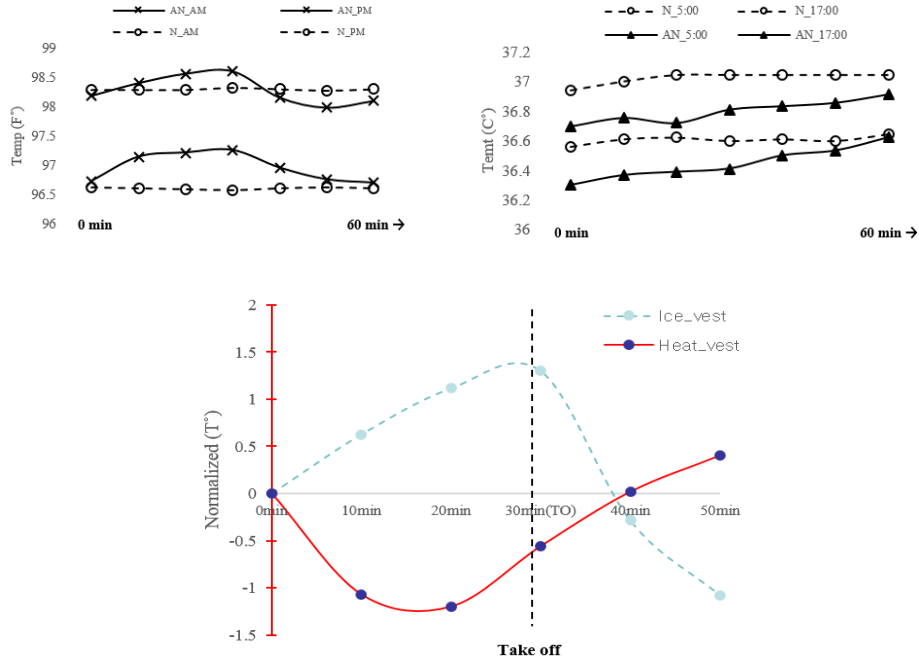

**Figure S14**. Illustration of the temperature according to the experimental design: normal (dotted line) and abnormal (dashed line) conditions. The horizontal axis denotes the temperature check time (Take off = without the ice/heat vest under the perturbation condition). The vertical axis is the level of the temperature change as calculated in Fahrenheit (F°) and Celsius (C°). Upper left = separate temperatures between the normal (N: a.m. and p.m.) and heat-based abnormal (AN: a.m. and p.m.) conditions; upper right = separate temperatures between the normal (N: a.m. and p.m.) and ice-based abnormal (AN: a.m. and p.m.) conditions. Note: in the abnormal session, the data were collected after 30 minutes (from 30 minutes to 60 minutes) from taking the ice (or heat) vest off. Bottom = normalized temperature adaptation tendency according to the artificially managed perturbation (the dotted line denotes the ice-vest perturbation; the bolded line denotes the heat-vest perturbation).

This temperature perturbation would be interpreted as a thermodynamic variable; that is, it is not itself restricted to the usual set of thermodynamic variables, such as the mean internal energy and entropy.

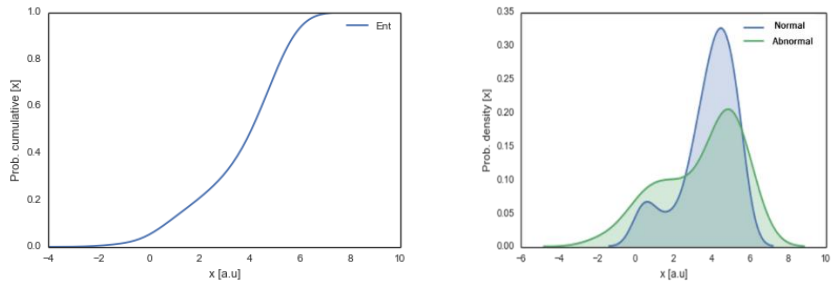

**Figure S15**. Estimated entropy production depending on the perturbation. The graph on the left side denotes the nominal distribution of all cases of entropy (x) production with a cumulative function (proportion [vertical axis] of the entropy value [horizontal axis]) with the arbitrary unit (a.u.). The graph on the right side represents the density of the entropy (data double plotted kernel

density function) according to the temperature perturbation (normal vs abnormal).

As shown in Figure S15, the emergence of the collective behavior of the increasing distance between $T(\chi, t) - T_0(t)$ via a perturbation will be related to an increase (or decrease) in the entropy. Comparing the states of $T_0 < \epsilon$ and $T_0 > \epsilon$, as shown in the graph on the right, as the distance-ordered state becomes smaller, the peak in the entropy production rate becomes higher. It is also important to note that the biological non-equilibrium bias toward a different temperature component hints at a possible deep connection between physical stability and entropy production. Given that the above entropy production is embedded in the order-disorder biophysical dynamic, the differences between the two circadian (termed a nearly 24-hour instance of oscillatory variation) rhythm points (am and pm), and the difference between the two temperature conditions in the psychomotor vigilances (i.e., $\dot{\phi}$ biological motor stability) were compared in order to determine whether differences in biological disorder resulted in differences in environmental perturbations.

According to the statistical testing of the data with direct reference to the research questions or hypotheses, the value of the equation "entropy = H (x)" was estimated considering the following null ($H_0$: $\theta_A = \theta_B$) and alternative ($H_a$: $\theta_A \neq \theta_B$) hypotheses. The hypothesis was proved that the different experimental conditions of the external source have a significant effect on the internal source (bi-manual motor variable) of $\dot{\phi}$. More specifically, different external components of the circadian processes or temperature have a significant effect on the internal stability; and the internal perturbation, from an external source, will have a significant effect on the degree of biological entropy. Explicitly, the statistic $F$ is calculated by dividing the difference between the group ($MS_{between}$) value by the difference between the value for the subjects within the group ($MS$). Observations are interested in the main effect of the circadian rhythm ($\alpha$), the temperature perturbation ($\beta$), and the interaction between the circadian rhythm and the temperature perturbation ($\alpha \times \beta$), as it affects the dependent variable of entropy production. Thus, the $F$ distribution was compared associated with each feature of interest to the error variance in order to determine if each effect is meaningful.

### 2.3. Results supplements

**Results for design 1:** Data collected from the participants (10: M=6, F=2, age 25 $\pm$ 3) at Seoul National University were used, but in order to measure the uncertainty, the calculation only involved the wrist joint data for entropy production. Each participant had to present four times with six trials each time (one participant, one minute, 24 trials = 6 [trials] *4 [circadian temperature points]).

**Table S5.** Each type of entropy value in normal day-night temperature effects.

| Participants (index) | Circadian 5:00 | Circadian 12:00 | Circadian 17:00 | Circadian 00:00 |
|---|---|---|---|---|
| $P1\_I_W$ | 3.92 | 2.88 | 4.02 | 5.82 |
| $P2\_I_W$ | 5.83 | 5.35 | 4.04 | 5.38 |
| $P3\_I_W$ | 4.19 | 3.82 | 3.69 | 4.73 |
| $P4\_I_W$ | 5.88 | 5.89 | 4.52 | 3.90 |
| $P5\_I_W$ | 5.78 | 5.25 | 4.43 | 4.31 |

| | | | |
|---|---|---|---|
| $P6\_$I$_W$ | 4.91 | 5.85 | 5.77 | 5.79 |
| $P7\_$I$_W$ | 5.76 | 5.62 | 4.06 | 5.83 |
| $P8\_$I$_W$ | 5.69 | 5.41 | 5.81 | 5.84 |

Note: *P* is the participant with each number of 1 ~ 8, W denotes the wrist dataset actually used from three different joint datasets, H (x) = entropy production. Note: the value of I is derived from the execution of trial ($w_1$, $w_2$), and the values from the two trials were divided by 2.

**Table S6.** Averaged entropy production from the normal day-night temperature values.

| | Circadian 5:00 | Circadian 12:00 | Circadian 17:00 | Circadian 00:00 |
|---|---|---|---|---|
| *N(I)* | 8 | 8 | 8 | 8 |
| *AVER(H)* | 5.246 | 5.010 | 4.544 | 5.200 |
| *STDEV(H)* | 0.796 | 1.078 | 0.812 | 0.782 |
| *SES* | 0.281 | 0.381 | 0.287 | 0.276 |

Note: *N(I)* = number of case indexed by the calculation of ($w_1$ + $w_2$ / 2), *AVER* = averaged entropy production; *STDEV* = averaged variability from the entropy production; *SES* = standard error score.

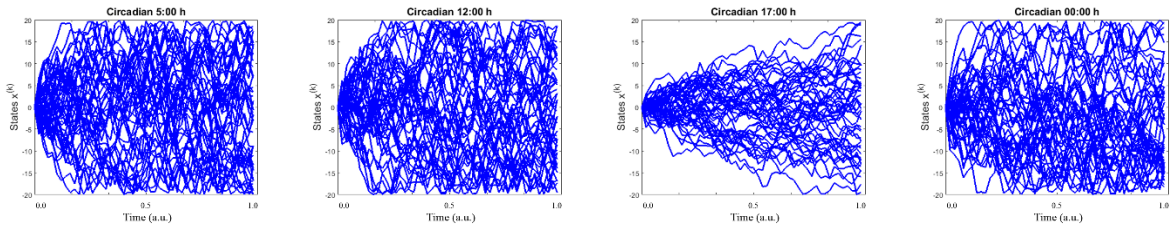

**Figure S16.** Uncertainty characteristics in the normal condition. The figures denote the estimated entropy states according to the time series. Left = 5:00, middle left = 12:00, middle right 17:00, right = 00:00. Note: for this realization, the data of participant 3 was used, which were representative of the scores most similar to the arbitrarily normalized scores gained by averaging all participants' scores as a function of the frequency competition.

The result shows the general feature of the average trend in ordinary circadian cycles. As shown in the table, the main effect of the uncertainty [H(x)] on the ongoing circadian cycle was not significant [$F(1, 3) = 1.074$, $\eta^2 = .823$, ($p < 0.376$)]. Absolute differences in the widths of the circadian cycle between temperature and entropy production can be observed, especially at the circadian points of 5:00 and 17:00 (t = 1.764, *p* < 0.103).

***Results for design 2:*** The data collected from the participants (10: M=6, F=2, age 25 $\pm$ 3) at Seoul National University were used, but only the wrist joint data were calculated with respect to entropy production. The data for each participant were analyzed four times, with six trials each time. Normal and abnormal day–night circadian temperature effects were compared at dawn (5 a.m.) and at dusk (5 p.m.) given that the core temperature reaches its maximum at around 5 p.m. and its minimum at approximately 5 a.m. In the perturbed condition, prior to the actual data collection, participants (n = 8) donned a heated vest for 30 minutes, which perturbed their core temperature. This additional data collection gave us the opportunity to compare the abnormal data to the previous normal data [one participant, one minute, 24 trials = 6 (trials) *2 (circadian points: normal data set) *2 (temperature perturbations dataset: body core temperature perturbed by the heated vest)].

**Table S7.** Each type of entropy value in normal and abnormal (heat-based) day–night temperature effects.

| Participants (index) | N_5:00 Normalized (Z) | N_17:00 Normalized (Z) | Ab_5:00 Normalized (Z) | Ab_17:00 Normalized (Z) |
|---|---|---|---|---|
| $P1$_$I_W$ | -0.67 | -0.59 | 0.88 | -0.99 |
| $P2$_$I_W$ | 0.89 | -0.58 | 0.9 | -0.91 |
| $P3$_$I_W$ | -0.45 | -0.86 | -0.01 | -0.59 |
| $P4$_$I_W$ | 0.93 | -0.18 | 0.87 | -2.33 |
| $P5$_$I_W$ | 0.84 | -0.26 | -0.78 | -0.51 |
| $P6$_$I_W$ | 0.14 | 0.84 | 0.91 | -1.09 |
| $P7$_$I_W$ | 0.83 | -0.56 | 0.89 | -0.2 |
| $P8$_$I_W$ | 0.77 | 0.87 | 0.83 | 0.16 |

Note: $P$ is the participant with each number of 1 ~ 8, W denotes the wrist dataset actually used from three different joint datasets. N means a normal circadian condition, and Ab denotes heat-vested abnormal circadian condition. Note: for more dramatic visualization, we used standard score (Z calculation). The value of $I$ is derived from the execution of each trial ($w_1$, $w_2$), and the value of these two trials' value was divided by 2.

**Table S8.** Averaged entropy production from normal and abnormal (heat-based) day–night temperature effects.

| | Circadian N_5:00 | Circadian N_17:00 | Circadian Ab_5:00 | Circadian Ab_17:00 |
|---|---|---|---|---|
| $N(I)$ | 8 | 8 | 8 | 8 |
| $AVER(H)$ | 0.410 | -0.165 | 0.564 | -0.809 |
| $STDEV(H)$ | 0.651 | 0.664 | 0.627 | 0.745 |
| $SES$ | 0.230 | 0.235 | 0.222 | 0.264 |

Note: $N(I)$ = number of case indexed by the calculation of ($w_1 + w_2$ / 2), $AVER$ = averaged entropy production; $STDEV$ = averaged variability from the entropy production; $SES$ = standard error score.

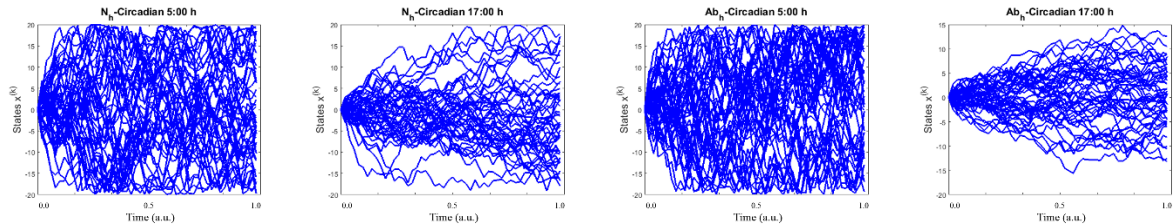

**Figure S17.** Uncertainty characteristics in the normal vs abnormal (heat-based) conditions. The figures denote the estimated entropy states according to the time series between normal (5:00 and 17:00) vs abnormal (5:00 and 17:00). N = normal, Ab = abnormal in terms of heat-based experimental design. Note: for this realization, participant 2's data was used which were representative of the most similar scores to the arbitrarily normalized scores gained by averaging all participants' scores as a function of the frequency competition.

The result shows the biological stability depending on the circadian time point, including the temperature (artificially perturbed body core temperature caused by the heated vest) perturbation. As shown in this figure, the main effect in the temperature perturbations was [$F(1, 3) = 1.301$, $\eta^2 = 0.961$, ($p < 0.258$)]. The main circadian effect was [$F(1, 3) = 20.531$, $\eta^2 = 15.166$, ($p < 0.001$)], and the significant temperature perturbation by the circadian cycle on the biological motor synchrony disorder was [$F(1, 3) = 3.453$, $\eta^2 = 2.551$, ($p < 0.068$)]. These results indicate that although the participants exhibited significantly greater levels of entropy in 5:00 a.m. conditions compared to the 17:00 p.m. conditions, in both the normal and the abnormal conditions (circadian effect), the temperature-associated disorder difference between a.m. and p.m. was intensified during artificially increased body core temperature (interaction effect).

***Results for design 3:*** Data collected from the participants (8: M=5, F=3, age 25 $\pm$ 3) at the University of Connecticut were used, but only the wrist joint data were calculated regarding entropy production. The data for each participant were analyzed four times with six trials each time. Normal and abnormal day–night circadian temperature effects were compared at dawn (5 a.m.) and at dusk (5 p.m.) considering that our core temperature reaches its maximum at around 5 p.m. and its minimum at approximately 5 a.m. In the perturbed condition, prior to the actual data collection, participants (n = 8) donned an ice vest for 30 minutes, which perturbed their core temperature. This additional data collection also gave us the opportunity to compare the abnormal data with the previous normal data [one participant, one minute, 24 trials = 6 (trials) *2 (circadian points: normal data set) *2 (temperature perturbations dataset: body core temperature decreased by the ice vest)].

**Table S9.** Each type of entropy value in normal and abnormal (ice based) day–night temperature effects.

| Participants (index) | N_5:00 Normalized (Z) | N_17:00 Normalized (Z) | Ab_5:00 Normalized (Z) | Ab_17:00 Normalized (Z) |
|---|---|---|---|---|
| *P1*_I$_W$ | 0.39 | -1.91 | -0.32 | -1.92 |
| *P2*_ I$_W$ | 0.73 | 0.51 | -0.01 | -1.85 |
| *P3*_ I$_W$ | 0.52 | -0.42 | 1.05 | -2.01 |
| *P4*_ I$_W$ | 0.98 | 0.86 | 0.93 | 0.31 |
| *P5*_ I$_W$ | 0.31 | -0.01 | 0.73 | -0.54 |
| *P6*_ I$_W$ | 0.26 | 0.44 | 1.06 | -0.03 |
| *P7*_ I$_W$ | -0.53 | -1.51 | 0.95 | -0.97 |
| *P8*_ I$_W$ | 0.57 | 0.64 | 0.47 | 0.28 |

Note: *P* is the participant with each number of 1 ~ 8, W denotes the wrist dataset actually used from three different joint datasets. N means normal circadian condition, and Ab denotes ice-vested abnormal circadian condition. Note: in order to gain a more dramatic visualization, we used standard score (Z calculation). The value of *I* is derived from the execution of each trial ($w_1$, $w_2$), and the value of these two trials was divided by 2.

**Table S10.** Averaged entropy production from normal and abnormal (ice based) day–night temperature effects.

| | Circadian N_5:00 | Circadian N_17:00 | Circadian Ab_5:00 | Circadian Ab_17:00 |
|---|---|---|---|---|
| *N(I)* | 8 | 8 | 8 | 8 |
| *AVER(H)* | 0.404 | -0.172 | 0.608 | -0.840 |
| *STDEV(H)* | 0.446 | 1.031 | 0.518 | 0.993 |
| *SES* | 0.158 | 0.365 | 0.183 | 0.351 |

Note: *N(I)* = number of case indexed by the calculation of ($w_1 + w_2$ / 2), *AVER* = averaged entropy production; *STDEV* = averaged variability from the entropy production; *SES* = standard error score.

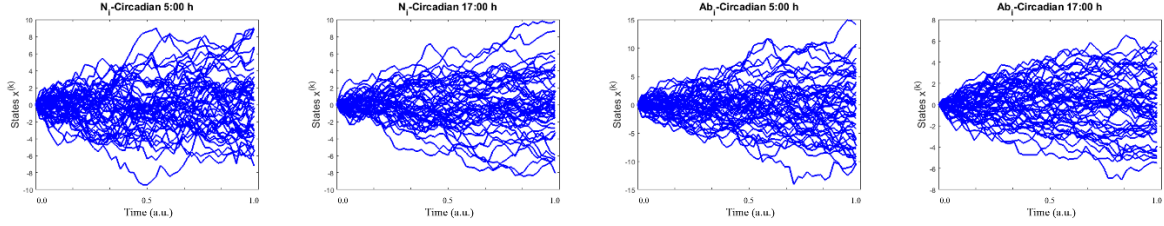

**Figure S18.** Uncertainty characteristics in the normal vs abnormal (ice based) conditions. The figures denote the estimated entropy states according to the time series between normal (5:00 and 17:00) vs abnormal (5:00 and 17:00). N = normal, Ab = abnormal in terms of ice-based experimental design. Note: for this realization, the data of participant 5 were used, as these were representative of the scores most similar to the arbitrarily normalized scores gained by averaging all participants' scores as a function of the frequency competition.

The result shows that biological stability is dependent on the circadian time point, including the temperature perturbation (artificially perturbed body core temperature cuased by the ice vest). As shown in the table, the main effect within the temperature perturbations was [$F(1, 3) = 1.211$, $\eta^2 = 0.861$, $(p < 0.275)$]. The circadian main effect was [$F(1, 3) = 23.041$, $\eta^2 = 43.317$, $(p < 0.001)$], and the significant temperature perturbation by the circadian cycle on the biological motor synchrony disorder was [$F(1, 3) = 4.264$, $\eta^2 = 3.035$, $(p < 0.043)$]. These results indicate that although the participants exhibited significantly higher levels of entropy in 5:00 a.m. conditions compared to the 17:00 p.m. conditions in both the normal and abnormal conditions (circadian effect), the temperature-associated disorder difference between a.m. and p.m. were intensified when the body core temperature was artificially perturbed (interaction effect).

### 2.4. Entropy calculation

Application of the average uncertainty model (*30*) to the actual data set of the experiment. Let us apply the above calculation to this experimental model. If we have a fair object $\phi$, of which the symbols have probabilities of 0 degrees = 0.5 and non-0 degrees = 0.5, we would calculate H of x as 0.5 multiplied by the log base 2 of 1 over 0.5, plus 0.5 multiplied by the log base 2 of 1 over 0.5, as follows:

$$\text{H}(\chi) = 0.5 \times log_2\left(\frac{1}{0.5}\right) + 0.5 \times log_2\left(\frac{1}{0.5}\right) \tag{8}$$

The result is 1. As

$$\frac{1}{0.5} = 2 \tag{9}$$

and

$$log_2(2) = 1 \tag{10}$$

then

$$H(\chi) = 0.5 \times 1 + 0.5 \times 1$$
$$\therefore H(\chi) = 1 \tag{11}$$

In this procedure, the fair object of $\phi$ has one bit of entropy, implying that when we calculate an object, we will receive an average of one bit of information. However, if the object is not fair, that is, when the symbols lead us to have probabilities of 0 degrees = 0.75 percent of the time but non-0 degrees = 0.25 percent of the time, the object for our information source will be

$$H(\chi) = 0.75 \times log_2\left(\frac{1}{0.75}\right) + 0.5 \times log_2\left(\frac{1}{0.25}\right) \tag{12}$$

and the result is approximately 0.811 because

$$\frac{1}{0.75} = 1.333 \ldots , \qquad \frac{1}{0.25} = 4 \tag{13}$$

and

$$log_2(1.333 \ldots) = 0.415, \qquad log_2(4) = 2 \tag{14}$$

then

$$H(\chi) = 0.75 \times 0.415 + 0.25 \times 2$$
$$\therefore H(\chi) = 0.811 \tag{15}$$

Thus, every time we measure an object of $\phi$ that is not fair, the calculation would give us an average of 0.811 bits of information, indicating that we will receive less (or more) of the information source (or entropy) from this un-faired object of $\phi$ than we would receive from a fair object of $\phi$.

Considering that the actual dataset will not be simply two states, if we extend the procedure, one can imagine an object's possibilities, as follows:

$$\{0.1, 0.1, 0.1, 0.5, 0.1, 0.1\}, \text{ or } \{0.16, 0.16, 0.16, 0.16, 0.16, 0.16\}$$

The same number of trials, but with different sets of both objects, will be approximated as follows:

Left set: $2.1 = 0.1 \times log_2\left(\frac{1}{0.1}\right) + 0.1 \times log_2\left(\frac{1}{0.1}\right) + 0.1 \times log_2\left(\frac{1}{0.1}\right) + 0.5 \times log_2\left(\frac{1}{0.5}\right) + 0.1 \times log_2\left(\frac{1}{0.1}\right) + 0.1 \times log_2\left(\frac{1}{0.1}\right)$

Right set: $2.5 = 0.16 \times log_2\left(\frac{1}{0.16}\right) + 0.16 \times log_2\left(\frac{1}{0.16}\right) + 0.16 \times log_2\left(\frac{1}{0.16}\right) + 0.16 \times log_2\left(\frac{1}{0.16}\right) + 0.16 \times log_2\left(\frac{1}{0.16}\right) + 0.16 \times log_2\left(\frac{1}{0.16}\right)$

This means that, on average, we will receive less (plot on the left = 2.1) entropy from a balanced object than from an unbalanced object (plot on the right = 2.5).

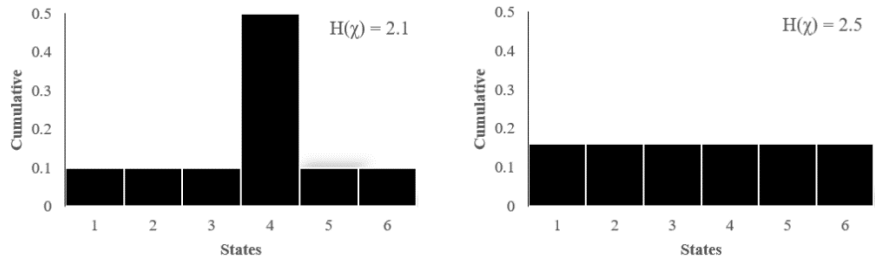

**Figure S19**. Simulation with different objects. The plot on the left denotes an un-faired set while the plot on the right denotes a faired set. Note: regarding the comparison of this simulation with variability, the variability of the left set, calculated by way of standard deviation, was 0.163, while the variability of the right set, also estimated by standard deviation, was 0.000. A comparison like this reflects that the uncertainty calculation shows us different information in spite of the same source being used.

Inspired by this simulation, we defined the probabilities as non-zero relative phase heights (cumulative function). This procedure allows us to calculate the uncertainty of an information source with any number of sets, noting that such a formulation used with the above procedures was based on observations made and experiments done on a macroscopic dimension (meaning any tangible piece of matter that we can see and work with within a laboratory).

### 2.5. Preliminary Pilot Test for the Experiment

For the preliminary pilot test, in order to find a relevant in-phase bi-manual synchrony variable, the collected data were calculated to compare different characteristics between one joint performance and several different joint performances. Although the wrist point is properly used, compared to other possible bi-manual pendulum areas such as the elbow or shoulder, this value is representative of the overall characteristics of a system that must be assessed. Moreover, repeating the assessment for only one position under several different conditions and several trials is likely to be associated with learning (or fatigue) effects.

**Table S11.** Data collection for pilot experiment:

| Group | Participants (N) | Body joint | Trials | Task/rest (min) |
|-------|-----------------|------------|--------|-----------------|
| G1    | 8               | Wrist      | 6      | 1 m / 5 m       |
|       |                 | Wrist      | 2      | 1 m / 5 m       |
| G2    | 8               | Elbow      | 2      | 1 m / 5 m       |
|       |                 | Shoulder   | 2      | 1 m / 5 m       |

Note: Group one = 8 participants, 6 trials at wrist point (total dataset = 48). Group 2 = 8 participants, 3 joints, and 2 trials at each joint with random sequences (total dataset = 48). Duration of each trial is 1 minute with a 5-minute rest interval between trials.

Based on these results, a means of data collection (excluding trial effects) was created as a relevant dependent variable that can be used to measure the internal source of stability. This implies a demonstration of the question as to whether the data, which were from only one joint in many trials, can be represented as the well-defined characteristic of the system.

The sessions were divided into two conditions (one joint = group 1, and different joint = group 2) with only one normal-temperature embedding cycle (5:00 pm: highest peak circadian rhythm of the core temperature over amplitude with skin capacitance). Each trial block lasted 1 minute with a 5-minute rest. During the one-joint session (wrist), participants received instruction about the preferred pendulum movements to establish in-phase 1:1 frequency locking at a 1.21 s metronome beat; this period was chosen because it corresponded to the natural period of the pendulum system without concern over amplitude or frequency (*29*). In the different joint sessions (wrist, elbow, and shoulder), the participants received instruction about the preferred pendulum, as in the single-joint session, but with the additional instruction of keeping different joints voluntarily fixed. A small amount of experience was provided with instruction beforehand to help avoid problems in complying with the session requirements. During the actual trial, no feedback was given, and the participants were asked not to report except when a problem arose. If the participants accidentally moved a joint that was supposed to be voluntarily fixed, the data from that trial were not analyzed, and the trial was repeated at a later time. The three different oscillation joints were used in random order.

Data were collected from 16 participants (at the University of Connecticut) (M=10, F=6, age 22 $\pm$ 3) in the normal condition (5:00 p.m.) in order to compare one typical anatomical position (wrist: M=5, F=3) and several different joint positions (wrist, elbow, and shoulder: M=5, F=3). Participants were divided into different experimental groups and asked to engage in the bimanual coordination in-phase 1:1 frequency locking at a 1.21 s metronome beat. This was for the following two reasons: (a) compared with other possible bi-manual pendulum areas, such as the elbow and the shoulder, the wrist point is commonly considered to represent this value, as the overall characteristic of the system must be confirmed. Moreover, (b) there was a doubt that repeating only one position under several different conditions through a number of trials could be associated with learning (or fatigue) effects.

First, investigation was divided into the phases for the different conditions in line with the requirements of the trial. Then, the data were analyzed to two values under the same normal condition (5:00 p.m.) with regard to which one would be the best experimental dependent variable to illustrate the different trial effects between the two conditions of the one joint (wrist) or the position of different joints (wrist, elbow, and shoulder).

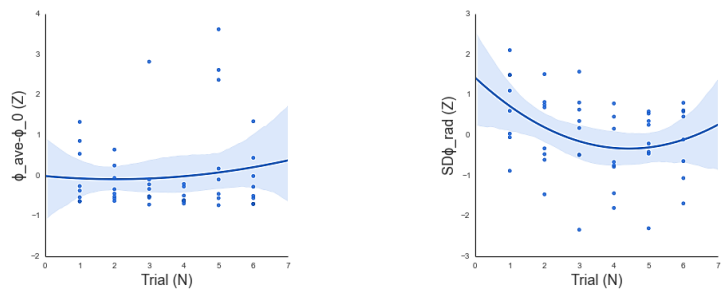

**Figure S20.1**. Repeated measure of the wrist according to the trial. Left = from the intended phase, $\phi_{ave} - \phi_0(rad)$, and right = standard deviation of the relative phase, $SD\phi(rad)$. Z is the standard score of the observed raw score $\chi$ (formula: $Z = \frac{\chi - \mu}{\sigma}$).

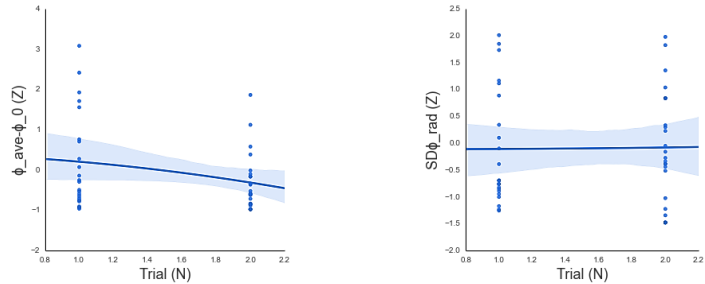

**Figure S20.2.** Repeated measure of the different joint (wrist, elbow, and shoulder) according to the trial. Left = from the intended phase, $\phi_{ave} - \phi_0(rad)$, and standard deviation of relative phase, $SD\phi(rad)$ = right. Z is the standard score of the observed raw score $\chi$ (formula: $Z = \frac{\chi - \mu}{\sigma}$).

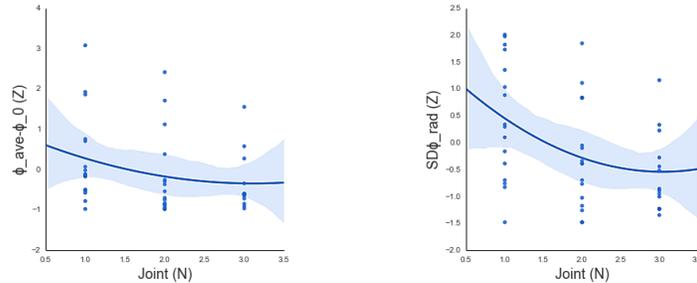

**Figure S20.3.** Deviation of the phase for the different joint orientation (1 = wrist, 2 = elbow, 3 = shoulder). Left denotes the topological effect [left = mean relative phase from the intended phase, $\phi_{ave} - \phi_0(rad)$, right = standard deviation of relative phase, $SD\phi(rad)$]. Z is the standard score of the observed raw score $\chi$ (formula: $Z = \frac{\chi - \mu}{\sigma}$). Note: In addition, we calculated the correlation taking into consideration the interaction between the two values and we found that both conditions have a significant relationship [$\phi_{ave} - \phi_0(rad)$] and [$SD\phi(rad)$]. This indicates that although the autocorrelation functions were different according to time series between [$\phi_{ave} - \phi_0(rad)$], [$\phi_{ave} - \phi_0(rad)$], [$SD\phi(rad)$], the higher [$\phi_{ave} - \phi_0(rad)$] correlated with higher [$SD\phi(rad)$] and lower [$\phi_{ave} - \phi_0(rad)$] also correlated with lower [$SD\phi(rad)$]: [Wrist Pearson Correlation $R = .46$ ($p = 0.0011$)], [Topology Pearson Correlation $R = .5$ ($p = 0.00041$)]. Such characteristics correspond to our predicted illustration of the relative phase based on the coordination dynamic calculations: $\phi_{ave} - \phi_0$ = fixed point shift, $SD\phi$ = variability as a function of frequency competition.

The pilot results showed that motor performance varied according to the different anatomical parts given the shift of mean relative phase from an intended phase [$\phi_{ave} - \phi_0(rad)$] and the in-phase variability [$SD\phi(rad)$] of these three behavioral variables. The important finding from these observations is that repetition with one joint may be significantly associated with decreased variability, akin to trial effects (see Figure 16.1) [Pearson Correlation $R = -.284$ ($p < 0.025$)], while the other condition (different topology) does not have a significance effect, according to the trial [Pearson Correlation $R = .110$ ($p < 0.236$)]. However, there is hierarchical significance, as the mean relative phase from an intended phase [$\phi_{ave} - \phi_0(rad)$] and the variability [$SD\phi(rad)$] were significantly wider for distal anatomy: $F(2, 47) = 4656.999$, $\eta^2 = 8.077$ ($p < 0.001$).

This investigation specifically considers these differences in order to determine the fundamental characteristic of different anatomical joints.

### Wrist       Elbow       Shoulder

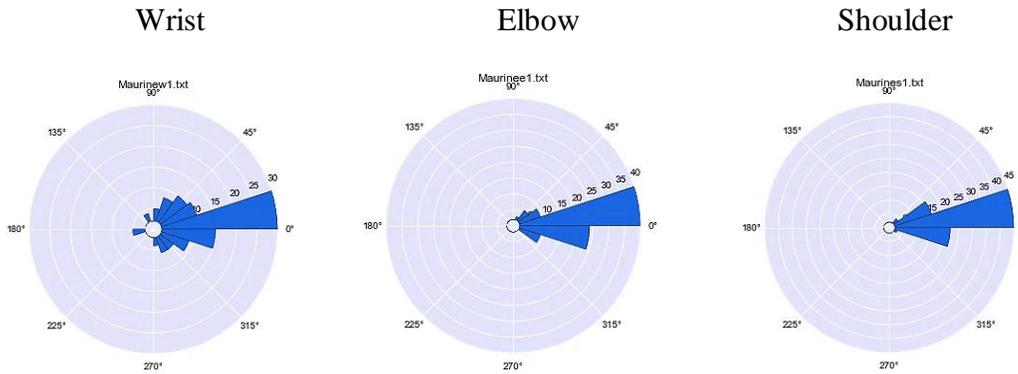

**Figure S21.** Circular representation of the different joint–coupled oscillations: Left = wrist, middle = elbow, right = shoulder (shaded sector arear means pendulum angle degree and variance). Note: with respect to these circular functions, we used 2π as a default (0) and calculated x using (180 degree*x/pi). In the above sample case, the degree represents how close in in-phase (0 degree, or 360 degree) or anti-phase (±180 degree) with variance (distribution).

Figure S21 illustrates the average stabilities of each joint for the wrist, elbow, and shoulder performances. Each topological asymmetry, couplings, and noise oscillations are reflected, showing that they were significantly wider for the distal (wrist) than the proximal (shoulder) with different anatomical parameters. This characteristic of symmetrical bimanual relationships may indicate greater heterogeneity of the scaling exponents at certain topological point.
Inspired by these analyses, it was decided to collect different values for three joints but to use only the wrist data to ascertain a biological characterization. This was for the following reasons: (a) although there was a significant learning effect when the participants undertook a task repeatedly with only the wrist joint, this led to significant typicality compared to the other two (elbow and shoulder) datasets. In order to overcome the learning (or fatigue) effects, the data were collected for three joints randomly, but only the wrist data were used. (b) Of course, although investigations were able to use combined data, which included all three different motor positions, representation of the combined data as the characteristic of a biological system appears to remove the important value of representativeness, as this combination likely has too many variables to manage. Moreover, (c) there was an expectation that collecting different motor scales but using widely represented data (wrist position) may meet both requirements of typicality as a well-defined system characteristic, as well as eliminating the learning effect stemming from the numerous trials. Thus, this manner was chosen to represent a typical internal source (dependent variable).

### *2.6. Graphical illustration of the inclusion and exclusion case of the data*

Graphical illustration of the inclusion and exclusion criteria (Upper = inclusion case, middle and bottom = exclusion case).

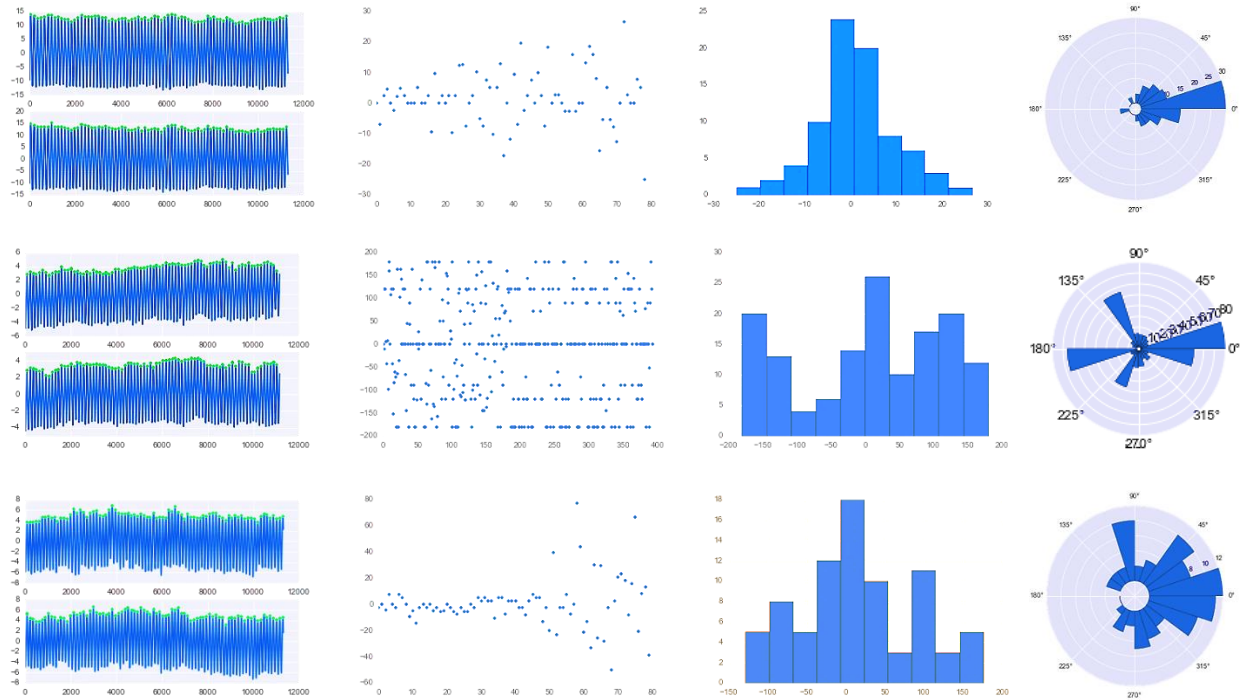

**Figure S22.** Graphical illustration of the inclusion and exclusion criteria. Note: to show the exclusion criteria, we added the data on exclusion cases based on the above procedure. As shown in the figure, even if same in-phase 1:1 frequency locking was applied, data from the task were likely contaminated by distractions from participants (figures at bottom) or devices (figures in the middle). In such a case, based on this basic rudimentary calculation, we excluded the data from five participants (Experiment 1 = 2, Experiment 2 and 3 = 3). The figure on the left denotes the frequency range of the amplitude (horizontal axis = time series and vertical axis = displacement, with the upper figure denoting the left-hand side and the bottom denoting the right-hand side). The figure on the left in the middle illustrates the discrete relative phase synchrony (horizontal axis = time series, vertical axis = relative phase checked peaks). The figure on the right in the middle is a phase histogram [horizontal axis = relative phase (0 equal in phase, $\pm$ 180 equal $\pm$ anti phase), vertical axis = proportion of the occurrence]. The shaded section in the figure on the right represents the pendulum angle degree and variance [with respect to this circular function, we used $2\pi$ as the default value (0) and calculated x via (180 degree*x/pi). In the sample case above, the degree represents the degree of closeness to the in-phase (0 degree, or 360 degree) or anti-phase ($\pm$180 degree) conditions and the distribution of the joints' relative phases].